\documentclass{article}

\usepackage{PRIMEarxiv}

\usepackage[utf8]{inputenc} 
\usepackage[T1]{fontenc}    
\usepackage{hyperref}       
\usepackage{url}            
\usepackage{booktabs}       
\usepackage{amsfonts}       
\usepackage{nicefrac}       
\usepackage{microtype}      
\usepackage{lipsum}
\usepackage{fancyhdr}       
\usepackage{graphicx}       
\usepackage{dsfont}
\usepackage{float}
\usepackage{tabularx}
\usepackage{subfig}
\usepackage{amsmath}
\usepackage{xcolor}
\usepackage[linesnumbered,lined,boxed,commentsnumbered]{algorithm2e}
\usepackage{epstopdf}
\newtheorem{remark}{Remark}%
\newtheorem{assumption}{Assumption}%
\newtheorem{conjecture}{Conjecture}%
\title{Particle Filter-Based Optimization: A Bayesian Approach for Global Stochastic Optimization
}

\author{
  Mostafa Eslami, Maryam Babazadeh \\
   Department of Electrical Engineering\\
  Sharif University of Technology \\
  Tehran, Iran\\
  \texttt{\{mostafa.eslami, babazadeh\}sharif@edu} \\
}

\begin{document}

\maketitle

\begin{abstract}
This paper proposes a novel global optimization algorithm, Particle Filter-Based Optimization (PFO), designed for a class of stochastic optimization problems in which the objective function lacks an analytical form and is subject to noisy evaluations. PFO utilizes the Bayesian inference framework of Particle Filters (PF) by reformulating the optimization task as a state estimation problem. In this context, evaluations of the objective function are interpreted as measurements, and a customized transition model based on covariance ellipsoids is introduced to guide particle propagation. This model serves as a surrogate for classical acquisition functions, equipping the PF framework with local search capabilities and supporting efficient exploration of the global optimum. To mitigate the adverse effects of measurement noise, the Unscented Transform (UT) is employed to approximate the underlying mean of the objective function, enhancing the accuracy of particle updates. The algorithm offers notable improvements over existing stochastic optimization algorithms for black-box multi-modal objective functions. First, PFO provides a fully probabilistic definition of particle weights, enhancing adaptability and robustness. Second, PFO integrates exploration and exploitation within a unified Bayesian framework, ensuring a non-zero probability of sampling from unexplored regions throughout the optimization process. This approach contrasts with traditional particle filter methods that are primarily used for state estimation, and heuristic optimization algorithms that lack theoretical guarantees. The novelty of PFO lies in its unique integration of particle filtering with a dynamic search space prediction, offering a theoretically grounded alternative to acquisition functions in Bayesian Optimization (BO).
\end{abstract}


\section{Introduction}\label{intro}
Global optimization is a subject of tremendous potential application, encompassing numerous fields such as engineering, economics, and artificial intelligence. Despite significant research and practical advancements over the past two decades, progress in the computational aspects of global optimization has not kept pace with improvements in digital computing power or the growing breadth of potential applications. This discrepancy is largely due to a persistent gap between theoretical developments and practical implementations, particularly between mathematical and heuristic approaches \cite{mockus2012bayesian}. 

The practical significance and inherent complexity of global optimization have motivated the development of a wide range of methodologies. These methods can generally be categorized into heuristic and non-heuristic approaches, although in practice this boundary is often blurred. Heuristic algorithms, which are designed to deliver acceptable solutions within feasible computational timeframes, have gained widespread attention. This popularity is attributed to the relative underdevelopment of rigorous mathematical theories for global optimization compared to local optimization \cite{zhigljavsky2007stochastic}. Consequently, global optimization continues to attract researchers across disciplines who seek to address the challenges posed by difficult real-world problems.

Efforts to bridge the gap between theory and practice increasingly emphasize integrating both theoretical analysis and empirical methodologies. This integration should not only focus on synthetic benchmark problems but also consider real-world examples under controlled assumptions and transparent formulations \cite{mockus2012bayesian}. In the context of noisy black-box functions, heuristic algorithms are generally expected to yield solutions sufficiently close to the global optimum rather than guaranteeing exact optimality \cite{simon2013evolutionary}.

Among heuristic strategies, stochastic optimization methods are of particular interest. Unlike deterministic algorithms, stochastic optimization incorporates randomness in various forms, such as noisy evaluations of the objective function, probabilistic decision rules, or random search directions \cite{zhigljavsky2007stochastic}. The present work focuses on stochastic optimization problems where objective function evaluations are corrupted by random errors, specifically addressing a class of NP problems denoted by $\mathcal{N}$. Extensive literature documents stochastic optimization algorithms for this subclass \cite{kingma2014adam,marti2005stochastic}.

Bayesian Optimization (BO) has emerged as a promising framework in this domain, utilizing prior and posterior distributions to iteratively guide the search toward global optima. BO is particularly effective for expensive black-box problems, whether noisy or noiseless, where evaluations are computationally costly \cite{chen2023hierarchical}. By modeling prior beliefs over the objective function and updating them based on observed data, BO selects new points using an acquisition function that balances exploration and exploitation \cite{chen2021function,mak2018efficient,dasgupta2008statistical}.

In parallel, estimation theory has contributed medium-term strategies in stochastic optimization by modeling objective function evaluations as noisy measurements. Algorithms such as the Heuristic Kalman Algorithm (HKA) and the Simulated Kalman Filter (SKF) utilize recursive state estimation concepts to locate optimal solutions \cite{toscano2009heuristic,ibrahim2015kalman}. Despite their practicality, these methods often lack convergence guarantees to global optima due to their heuristic nature and limitations in handling complex, multi-modal landscapes \cite{yang2010nature}.

An alternative line of research reformulates global optimization as an optimal control problem, where particle-based methodologies, such as the Feedback Particle Filter (FPF), are applied to govern the evolution of candidate solutions \cite{zhang2018mean,zhang2017controlledthesis}. These approaches establish theoretical connections between particle filtering and mean-field control, revealing that density transport in the state space can be interpreted as a gradient flow of an optimal value function. While this provides a robust theoretical foundation and enables advantages such as bypassing resampling, the method proposed here differs by directly employing particle filters with a novel transition model based on covariance ellipsoids. This alternative focuses on direct integration of exploration and exploitation and introduces a new mechanism to address noise using the Unscented Transform (UT).

While several methods have adopted particle filtering for stochastic global optimization \cite{stinis2012stochastic,zhou2013particle,gerber2022global}, they differ significantly in core methodology. These methods typically manipulate aspects of the PF algorithm to enhance convergence or unify various randomized optimization techniques through the lens of filtering theory. They often rely on modifying control parameter assignments for particle offspring. In contrast, the proposed method explicitly addresses noisy, nonlinear objective evaluations by employing the UT to approximate the mean and covariance at selected points. This improves the accuracy of particle updates and introduces a novel transition model based on local covariance ellipsoids. These ellipsoids function similarly to acquisition functions in BO, maintaining the balance between exploration and exploitation and mitigating degeneracy caused by unknown transition dynamics. The method also accommodates alternative dynamics through selective updates of the ellipsoidal structure, whereas prior work typically incorporates random walks or uses the objective function directly within the state transition dynamics \cite{stinis2012stochastic,zhou2013particle}.

Other approaches designed for the same class of optimization problems, such as the parallel Sequential Monte Carlo (SMC) optimizer \cite{akyildiz2020parallel}, emphasize convergence guarantees and theoretical rates. Although based on a similar particle-driven concept, these methods primarily focus on the SMC framework and can be adapted into the proposed Particle Filter-Based Optimization (PFO) structure as a replacement for the PF component.

Recent studies in machine learning optimization, especially in large-scale empirical risk minimization, have also explored PF-based algorithms. For instance, Particle Filter-based Stochastic Optimizers (PFSOs) have been proposed to address the computational infeasibility of gradient-based methods on massive datasets. These methods integrate PF with the Incremental Proximity Method (IPM), although their performance deteriorates in highly nonlinear problems due to the reliance on numerical solvers \cite{liu2020particle,bertsekas2011incremental}.

The method introduced in this paper proposes a new stochastic global optimization algorithm that combines prior optimization states with posterior distributions of objective values. It features a predictive mechanism for estimating promising search regions and a theoretically grounded likelihood function for updating particle positions. A key innovation lies in the replacement of the traditional acquisition function from BO with a dynamic system perspective. The optimization process is guided by the natural gradient-like behavior of the objective function, modulated through a utility function applied over sigma points generated via the UT. As output measurements are noisy, UT is used to estimate the true mean more accurately \cite{simon2006optimal}. This also mitigates the impact of measurement errors during particle propagation and weight assignment.

Compared to Principal Component Analysis (PCA)-based stochastic optimization \cite{kuznetsova2012pca}, the UT provides a more statistically optimal approximation of the local mean and covariance. The particle update direction and magnitude are derived from the eigenstructure of the covariance matrix associated with each particle’s sigma points.

This probabilistic, locally adaptive search framework allows the proposed algorithm to explore both promising and non-promising regions effectively. Particles are treated as agents within a filter-theoretic model that preserves the ability to explore unvisited spaces while converging toward global optima. In this context, particle weights are interpreted through the posterior likelihood based on the deviation from the best-known objective value. The resulting method, termed Particle Filter-Based Optimization (PFO), offers a theoretically grounded alternative that integrates estimation theory with global optimization.

PFO is applicable to a broad class of black-box stochastic optimization problems, including machine learning applications such as empirical risk minimization with non-convex losses \cite{liu2020particle}, simulation-based optimization in engineering, finance, and logistics, and control systems tuning and motion tracking problems \cite{zhang2018mean}.

This paper is organized as follows. Section~\ref{sec:stoch} formulates the stochastic optimization problem, outlines the underlying assumptions, and revisits Particle Swarm Optimization (PSO) \cite{kennedy2010particle} from a filtering perspective. Section~\ref{sec:pfo} presents the proposed PFO framework with emphasis on the novel transition and update mechanisms. Section~\ref{sec:eval} evaluates the performance of PFO on standard benchmarks and includes a sensitivity analysis with randomized parameters, with results compared against PSO.

{\renewcommand\arraystretch{1}
\begin{table}[!h]
\begin{tabular}{p{.15\linewidth} p{.75\linewidth}}
\large{\textbf{Nomenclature}}&\\
&\\
$\mathds{R}^n$ &  Real values of dimension $n$ \\
$x$ & Optimization variable\\
$\mathcal{N}$ & Class of general stochastic optimization problems with randomness in evaluation of objective function $h(x)$\\
$\mathcal{C}$ & Sub-class of $\mathcal{N}$ with additive noise $v(x)$ to $h(x)$\\
$H(x)$ & Black-boxed objective function and member of $\mathcal{C}$, i.e. $H(x)=h(x)+v(x)$\\
$x_k$ & Evaluated optimization variable at iteration $k$-th\\
$\hat{x}_k$ & Estimation of $x$ at iteration $k$-th\\
$y_k$ & Measure of objective function  at iteration $k$-th, i.e. $y_k=H(x_k)$\\
$\hat{x}_k$ & Estimation of measurement objective function at iteration $k$-th\\
$p(A|B)$ & Conditional power density function of $A$ given $B$, the same definition is true for $q(A|B)$\\
$w_k$ & Weight of particles at iteration $k$-th\\
$(\cdot)^i$ & Superscript $i$ refers to $i$-th particle\\
$P^{xx}_k$ & Covariance matrix of random variable $x$ at iteration $k$-th\\
$C_{xy}$ & Cross covariance matrix of random variables $x$ and $y$\\
$\underline{\lambda}(A)$ & Minimum eigenvalue of matrix $A$\\
$\bar{\lambda}(A)$ & Maximum eigenvalue of matrix $A$\\
$\mathcal{N}(a,Q)$ & Normal distribution of mean $a$ and covariance $Q$\\
$\mathcal{X}^j\in \mathds{R}^{1\times n}$ & Sigma point $j$ of random variable $x$ in UT\\
$\mathcal{Y}^j\in \mathds{R}$ & Transformed sigma point $\mathcal{X}^j$\\
\end{tabular}
\end{table}
}
\newpage
\section{Problem Statement}\label{sec:stoch}
Stochastic optimization becomes essential when either the evaluation of the objective function \( h(x) \) is corrupted by random noise, or when the optimization algorithm itself incorporates randomness in the selection of search directions during the iterative solution process \cite{spall2005introduction}. 

Let \( \mathcal{N} \) denote the class of general stochastic optimization problems characterized by inherent randomness in the evaluation of the objective function. Let \( \mathcal{C} \subseteq \mathcal{N} \) represent the subclass of problems considered in this work. Let the feasible set for the optimization variable \( x \) be denoted by \( \mathds{D} \subseteq \mathds{R}^n \). The goal in problem class \( \mathcal{C} \) is to find the value(s) of \( x \in \mathds{D} \) that minimize a scalar-valued, noisy objective function \( h(x) \), whose evaluations are only accessible through corrupted measurements.

Let \( \hat{x}_k \) denote the estimate of the optimizer at iteration \( k \). Owing to the stochastic nature of the problem, \( \hat{x}_k \) is treated as a random vector. The noisy measurement model used throughout this paper is defined as,
\begin{align}
H(x) = h(x) + v(x),
\end{align}
where \( h(x) \) is the true, unknown objective function, and \( v(x) \) is a zero-mean, potentially state-dependent noise process that captures the uncertainty in the measurement. 
Let \( y_k \in \mathbb{R} \) denote the observed value of $H(x_k)$ at the $k$-th iteration, i.e., $y_k = H(x_k) = h(x_k) + v(x_k)$. Similarly, the objective value at the estimated optimization variable is denoted by $\hat{y}_k$, i.e.,  $\hat{y}_k = H(\hat{x}_k) = h(\hat{x}_k) + v(\hat{x}_k)$
\begin{assumption}
The analytic expression of $H(x)$ is unknown. However, point-wise evaluations of $H(x)$ can be obtained through noisy measurements.
\end{assumption}
\begin{assumption}
    The measurement noise is modeled as a zero-mean Gaussian process with possibly state-dependent covariance. Specifically, $v(x) \sim \mathcal{N}(0, R(x)),$
where $R(x)$ is a known, positive semi-definite covariance function that captures the possibly heteroscedastic nature of the measurement uncertainty.
\end{assumption}
The problem considered in this work belongs to a subclass of NP problems that are computationally expensive to solve but whose solutions can be verified efficiently. In some contexts, this class of problems is referred to as NP-complete \cite{zhigljavsky2007stochastic}.

A heuristic algorithm is designed to solve problems more quickly and efficiently than traditional methods by sacrificing optimality, accuracy, precision, or completeness for speed. Heuristic algorithms are often used to solve NP-complete problems, a class of decision problems where no known efficient way to find a solution quickly and accurately exists, although solutions can be verified when given. Heuristics can produce a solution individually or be used to provide a good baseline and be supplemented with optimization algorithms. They are most often employed when approximate solutions are sufficient and exact solutions are computationally expensive \cite{cook1983overview}. Such algorithms can find very good results without any guarantee of reaching the global optimum; often, there is no other choice but to use them.

There are generally two phases in solving NP problems using heuristic algorithms:
\begin{itemize}
\item Phase I (Diversification): This phase involves global exploration of the search space. The algorithm explores the entire domain to determine potentially good subregions for future investigation.
\item Phase II (Intensification): In this phase, local exploitation is performed.  Local optimization algorithms are applied to determine the final solution.
\end{itemize}

To illustrate the core concept, consider Algorithm \ref{alg:pso} which implements the PSO algorithm \cite{simon2013evolutionary}. The algorithm is based on the observation that groups of individuals work together to improve not only their collective performance on some tasks but also each individual's performance. It propagates its optimization variables' particles based on probabilistic velocity updates. Velocity is likely to change based on the best individual, neighbors, or global experiences. The level of influence is parameterized and determines the balance between diversification and intensification (lines 8 to 14). Following the velocity update, the new position of the particles updates to the minimum solution found so far for each. Finally, in line 16, the algorithm identifies the global best solution as the position associated with the lowest objective function value found among all particles.

The standard PSO algorithm \ref{alg:pso} is generally ineffective in solving problems in class \( \mathcal{C} \), as it lack the ability to distinguish between noise in the measurements and the true underlying objective values. Although PSO incorporates probabilistic velocity updates, the objective functions being minimized—whether local or global—are inherently deterministic. This mismatch limits the algorithm's performance in noisy settings. In \cite{taghiyeh2016new} this limitation is addressed by extending PSO frameworks to better handle noisy objective functions.

\begin{algorithm}[ht]
 \SetAlgoLined
 \KwData{$N$ [Number of particles], $\sigma$ [Number of nearest neighbors], 
$x_{\text{min}}, x_{\text{max}}$ [Position bounds], 
$v_{\text{min}}, v_{\text{max}}$ [Velocity bounds]}
 \KwResult{$g$ [global best solution], $y_b$ [function value at best solution]}
 \textbf{Initialize:} \\
\For{$i = 1:N$}{ Sample initial position: $x_i \sim \mathcal{U}(x_{\min}, x_{\max})$\;
Sample initial velocity: $v_i \sim \mathcal{U}(v_{\min}, v_{\max})$\;
Initialize personal best: $b_i \leftarrow x_i$\; }

Set global best: $g = \arg\min\{h(b_1), h(b_2), \ldots, h(b_N)\}$\;
 \While{termination criteria not met}{
 \For{$i = 1 : N$}{
 $H_i = \sigma$ nearest neighbors of $x_i$\;
 $h_i = \arg\min\{h(x): x\in H_i\}$\;
 Generate random vectors $\phi_p$, $\phi_n$ and $\phi_g$\;
 $v_i= v_i +\phi_p . (b_i-x_i) + \phi_n . (h_i-x_i) +\phi_g(g-x_i)$\;
 $x_i\leftarrow x_i+v_i$\;
 $b_i= \arg\min\{h(x_i),h(b_i)\}$\;
 }
 $g = \arg\min\{h(b_1), h(b_2), \ldots, h(b_N)\}$\;
 $y_b = h(g)$
 }
 \caption{Particle Swarm Optimization Algorithm: a Particle-Based Optimization Heuristic }\label{alg:pso}
\end{algorithm}

In addition, to effectively solve an optimization problem empirically, the balance between exploration and exploitation should be dynamically adjusted throughout the optimization process. Typically, a high degree of diversification is favored in the early stages to explore the search space broadly, while intensification becomes more prominent in later stages to refine promising solutions. Consequently, the influence of the best-known solutions must be adapted to the characteristics of the problem. However, improper tuning of this balance may lead to algorithmic instability or premature convergence to suboptimal local minima. For example, in Algorithm \ref{alg:pso} the influence factors, namely, $\phi_b$, $\phi_n$, and $\phi_g$ are randomly selected, without any statistical relationship between them or the velocities (next solutions or particle positions). 

In fact, the velocity update mechanism in PSO is analogous to the transition prior in a Bootstrap Particle Filter (BPF), elaborated in Section \ref{sec:pfo},  while the selection of the best individual shares similarities with the particle weight updates based on the measurement likelihood probabilities. Building on this conceptual alignment, Section \ref{sec:pfo} introduces a novel stochastic optimization algorithm, i.e., PFO which reinterprets the likelihood function from an optimization standpoint rather than a purely estimation-oriented view. PFO offers notable improvements over heuristic algorithms such as PSO in two key aspects. First, unlike PSO, all influence factors and weight assignments in PFO are probabilistically defined, enhancing both its adaptability and robustness. Second, PFO integrates exploration (diversification) and exploitation (intensification) within a unified framework, ensuring a non-zero probability of sampling from unexplored regions of the search space throughout the optimization process.

\section{PFO: Particle Filter-Based Optimization}\label{sec:pfo}
\subsection{Bootstrap Particle Filters}
Importance sampling is a general Monte-Carlo integration method that provides a recursive solution to nonlinear filtering problems using a Bayesian approach. The key idea in the PF is to represent the required posterior density function by a set of random samples with associated weights. Then, estimates are computed using these samples and weights. Samples evolve based on a proposal density function $q(x_k|x_{k-1},y_k)$, and weight updates are based on the equation,
\begin{align}\label{eq
}
w^i_k \propto \cdot w^i_{k-1} \dfrac{p(y_k|x_k^i)p(x_k^i|x_{k-1}^i)}{q(x_k^i|x_{k-1}^i,y_k)}.
\end{align}
The choice of importance density is one of the most critical issues in the design of a particle filter. The optimal importance density function that minimizes the variance of importance weights is $p(x_k|x_{k-1}^i,y_k)$ \cite{doucet2000sequential}. This posterior can be written for particle $i$ as,
\begin{align}\label{eq:BFP1}
p(x^i_k|x_{k-1}^i,y_k) = \dfrac{p(y_k|x^i_k,x_{k-1}^i)p(x^i_k|x_{k-1}^i)}{p(y_k|x_{k-1}^i)}.
\end{align}
Substituting \eqref{eq:BFP1} into \eqref{eq
} yields,
\begin{align}
w_k^i \propto w_{k-1}^i \cdot p(y_k|x_{k-1}^i).
\end{align}
These series of equations utilize sampling from the optimal proposal density and $p(y_k|x_{k-1}^i)$, which requires their analytical expressions. The analytical evolution of these posteriors is difficult in most cases, except for some special Gaussian problems \cite{doucet2000sequential}. Therefore, suboptimal methods that approximate the optimal importance density have been developed. The most popular choice is the transitional prior for the proposal density, i.e., $p(x_k|x_{k-1}^i)$. Substituting this into \eqref{eq
} yields:
\begin{align}
w_k^i \propto w_{k-1}^i \cdot p(y_k|x_k^i).
\end{align}
This type of PF is known as the Bootstrap Particle Filter (BPF), also known as the Sequential Importance Resampling (SIR) filter. BPF suffers from the lack of measurement in the transitional prior, which leads to the generation of unnecessary particles that are not in interest of the likelihood distribution. In the initial iterations, only a few particles will be assigned a high weight, causing particle degeneration. The Bootstrap Particle Filter assumes that,

\begin{assumption}
The state dynamics and measurement functions are assumed to be known.
\end{assumption}
\begin{assumption}
It is required to be able to sample realizations from the process noise distribution and from the prior.
\end{assumption}
\begin{assumption}
The likelihood function needs to be available for point-wise evaluation (at least up to proportionality).
\end{assumption}
A generic algorithm for BPF is presented in Algorithm \ref{alg:bpf}. These weak assumptions and the suboptimal choice of proposal density have encouraged researchers to slightly manipulate the generic procedure to obtain more efficient filters as variants of BPF. Methods like resampling, roughening, and regularizing have been developed for practical applications \cite{gordon2004beyond}. Methods exist to encourage the particles to be in the right place (in the region of high likelihood) by incorporating the current observation. One such method is the Auxiliary Particle Filter (ASIR), which introduces intermediate distributions between the prior and likelihood \cite{pitt1999filtering}. The basic idea in ASIR is to perform the resampling step at time $k-1$ (using the available measurement at time $k$), before the particles are propagated to time $k$. In this way, the ASIR filter attempts to mimic the sequence of steps carried out when the optimal importance density is available \cite{gordon2004beyond}.

The allowance for mimicking the optimal importance density via current measurement, along with the aforementioned assumptions, are the building blocks of the proposed stochastic optimization algorithm in Section \ref{sec:pfo_algorithm}.\\

\begin{algorithm}[ht]
 \SetAlgoLined
 \KwData{$N$ [Number of particles], $p(x_k|x^i_{k-1})$ [State transition model] and $p(y_k|x_k^i)$ [Observation likelihood]}
 \KwResult{$\hat{x}_k$ [State estimate] and $P_k^{xx} $  [Empirical covariance matrix]}
 \textbf{Initialization:} \\
\For{$i = 1 : N$}{
 Draw $x_0^i \sim \mathcal{U}(x_{\min}, x_{\max})$\;
    Assign $w_0^i = \frac{1}{N}$\;
}
\For{$k = 1, \dots, k_{max}$}{
 \For{$i = 1 : N$}{
 Draw $x_k^i\sim p(x_k|x^i_{k-1})$ \;
 Calculate $\tilde{w}^i_k = w_{k-1}^i p(y_k|x_k^i)$\;
 }
 $w_k^i = \tilde{w}^i_k/\sum_{i=1}^{N}\tilde{w}^i_k$\;
 
 $\hat{x}_k = \sum_{i=1}^{N}w_k^ix_k^i$\tcc*[r]{this is a MMSE estimation}
 
 $P_k^{xx} = \sum_{i=1}^{N} w_k^i(x_k^i-\hat{x}_k)(x_k^i-\hat{x}_k)^T$\tcc*[r]{this is empirical covariance matrix} 
}
\caption{Bootstrap Particle Filter}\label{alg:bpf}
\end{algorithm}

\subsection{Algorithmic Description of PFO}\label{sec:pfo_algorithm}
An alternative intuition for the PF presented here demonstrates a close connection between estimation and optimization problems. The PF propagates the last found representative of states ($x_{k-1}^i$ - optimization variables) through a transitional prior, i.e., $p(x_k^i|x_{k-1}^i)$, then harvests the most promising estimates by comparing the actual measurement ($y_k$) and estimated output ($\hat{y}_k$) based on the likelihood posterior. Since the states of the system under estimation are more likely to follow their dynamical trajectory, propagation through the transitional prior makes sense. However, if we could propagate based on the knowledge of the last obtained measurements, generating more particles close to the true solution could enhance estimation performance.

Consider an optimization problem where the agents (analogous to particles in a Particle Filter) are positioned at \( x_{k-1}^i \) at iteration \( k-1 \).  Similar to the estimation setting, the goal is to determine the best positions \( x_k^i \) by using both current and historical observations of the objective function. From an optimization perspective, this involves assigning each agent a trajectory that guides it from its initial random position toward the global minimum. This trajectory serves as the counterpart of system dynamics or the transition prior in the estimation framework.

By reinterpreting the concept of likelihood in the context of optimization, the PF algorithm can be transformed into a global optimization algorithm. Specifically, the likelihood is redefined as the probability associated with the error between the observed and predicted objective values, i.e.,
($y_k - \hat{y}_k$). Consequently, particle weights are assigned based on their proximity to the observed objective value. Over successive iterations, this probabilistic weighting, combined with local propagation guided by the transition prior, encourages convergence toward the global minimum. The resulting approach is referred to as the Particle Filter-Based Optimization algorithm, or PFO.

PFO serves as a global optimizer by uniquely combining diversification and intensification within a probabilistic framework. As discussed in Section \ref{sec:stoch}, for example, PSO evolves particles using velocity updates influenced by the best individual and collective experiences, followed by a deterministic selection of the best candidates. Therefore, the influences are not adapted to the past and current observations. Unlike PSO, the degree of influence and weight assignment in PFO are all probabilistic, enhancing adaptability and robustness. Also, the diversification phase in PFO is coupled with intensification and guarantees a non-zero probability of searching unexplored domains of the search space. In PFO, the most highly weighted particles predominantly contribute to estimating the global minimum. Specifically, particles associated with higher uncertainty explore the search space in pursuit of potential new minima, whereas particles with lower uncertainty receive greater weights and focus on exploiting the current estimated solution to enhance its accuracy. This mechanism effectively integrates the exploration–exploitation balance characteristic of particle filters, which is inherited in the proposed PFO Algorithm.

The proposed global optimization Algorithm PFO to solve  class-$\mathcal{C}$ problems is outlined as Algorithm \ref{alg:pfo}. At the initial iteration, a set of \( N \) particles \( \{{x}_0^i\}_{i=1}^N \) is drawn from a uniform prior distribution over the feasible domain, i.e., ${x}_0^i \sim \mathcal{U}({x}_{\min}, {x}_{\max})$. All particles are initialized with equal weights \( w_0^i = \frac{1}{N} \).\\
At each iteration \( k \), the algorithm begins by propagating each particle according to a stochastic state transition model. Specifically, for each particle \( {x}_{k-1}^i \), a new particle is sampled as
\({x}_k^i \sim p({x}_k \mid {x}_{k-1}^i)\).
Details on the selection of transitional prior is given in Section \ref{subsec:transitional}. The empirical local covariance is estimated as,
\[
{P}_k^{{x}^i} = ({x}_k^i - \hat{{x}}_{k-1})({x}_k^i - \hat{{x}}_{k-1})^\top + {Q},
\]
where \( \hat{{x}}_{k-1} \) is the weighted mean from the previous iteration and \( {Q} \) is a predefined exploration noise matrix. This local covariance is then used to perform an UT centered at \( {x}_k^i \) to approximate the transformed statistics of \( {x}_k^i \) under the nonlinear and potentially multi-modal objective function \( h(\cdot) \).

UT requires deterministic construction of \( 2n_x + 1 \) sigma points to describe the true mean and covariance of $x_k^i$. These sigma points and their corresponding weights are generated as,
\begin{align*}
\begin{aligned}
  &\mathcal{X}^{(0)} = x_k^i, \\
  &\mathcal{X}^{(j)} = x_k^i + \left(\sqrt{(n_x + \lambda)P_k^{x^i}}\right)_j,  \\
  &\mathcal{X}^{(j + n_x)} = x_k^i - \left(\sqrt{(n_x + \lambda)P_k^{x^i}}\right)_j,
\end{aligned}
\qquad
\begin{aligned}
  &w_0 = \frac{\lambda}{n_x + \lambda}, \\
  &w_j = \frac{1}{2(n_x + \lambda)}, \quad && j = 1, \ldots, n_x, \\
  &w_{j+n_x} = \frac{1}{2(n_x + \lambda)}, \quad && j = 1, \ldots, n_x.
\end{aligned}
\end{align*}
where \( \lambda  \) is a scaling parameter of the UT, and $\left(\sqrt{(n_x + \lambda)P_k^{x^i}}\right)_j,$ is the $j$th row of the matrix square root $L$ of $(n_x + \lambda)P_k^{x^i}$ such that $(n_x + \lambda)P_k^{x^i}=L^T L$. Each sigma point is passed through the objective function to yield transformed sigma values,
\[
\mathcal{Y}^{(j)} = h(\mathcal{X}^{(j)}), \quad j = 0, \ldots, 2n_x.
\]

Then, the first two moments of the output are computed as,
\[
{y}_k^i = \sum_{j=0}^{2n_x} w_j^{(m)} \mathcal{Y}^{(j)},
\]
\[
{P}_k^{{y}^i} = \sum_{j=0}^{2n_x} w_j^{(c)} \left(\mathcal{Y}^{(j)} - {y}_k^i\right)\left(\mathcal{Y}^{(j)} - {y}_k^i\right)^\top + {R_k},
\]
where \( R_k \) is the observation noise covariance for $v(x_k^i)$.
\begin{remark}
The spread of the sigma points from the mean value increases with the dimension \( n_x \). This spread can be adjusted to some extent by the choice of the scaling parameter $\lambda$. Moreover, the scaled UT is introduced to provide an additional degree of freedom in scaling the sigma points either away from or towards the mean values \cite{julier2002scaled}.
\end{remark}

Each particle is reweighted according to its likelihood under the estimated observation model,
\[
\tilde{w}_k^i = w_{k-1}^i \cdot p({y}_k \mid {x}_k^i),
\]
and normalized to produce posterior weights,
\[ w_k^i = \frac{\tilde{w}_k^i}{\sum_{j=1}^N \tilde{w}_k^j}. \]
Informally, particles whose states are closer to the most recent best estimate are likely to receive higher weights, whereas those farther away are assigned lower weights. The high-weight particles contribute more significantly to the next estimate of the optimal solution, while low-weight particles promote exploration of the search space for potential alternative minima.

\begin{algorithm}[p]\caption{Particle Filter-Based Optimization (PFO) with Unscented Transform}\label{alg:pfo}
\SetAlgoLined
\KwData{
$N$ [Number of particles], 
$p(x_k | x_{k-1}^i)$ [State transition model],  
$P_{\text{min}}^{xx}, P_{\text{min}}^{yy}$ [Minimum covariance thresholds], 
$k_{\text{max}}$ [Maximum iterations], $\lambda$ [Scaling parameter], $Q$ [Exploration noise matrix]
}
\KwResult{$\hat{x}_k$, $\hat{y}_k$ }

\textbf{Initialization:} 

\For{$i = 1 : N$}{
    Draw $x_0^i \sim \mathcal{U}(x_{\min}, x_{\max})$ \tcc*[r]{Uniform initialization}
    Set $w_0^i = \frac{1}{N}$
}

\For{$k = 1, 2, \ldots, k_{\text{max}}$}{
    \For{$i = 1 : N$}{
        Sample $x_k^i \sim p(x_k | x_{k-1}^i)$ \tcc*[r]{Details on transitional prior is given in Section \ref{subsec:transitional}}
        
        $P_k^{x^i} = (x_k^i - \hat{x}_{k-1})(x_k^i - \hat{x}_{k-1})^T + Q$ 
        
        $\mathcal{X}^{(0)} = x_k^i$ \tcc*[r]{Generate $2n_x + 1$ sigma points} 
        
        \For{$j = 1 : n_x$ }{
            $\mathcal{X}^{(j)} = x_k^i +  \left(\sqrt{(n_x + \lambda)P_k^{x^i}}\right)_j$  \tcc*[r]{\hspace{-2mm}$\left(A\right)_j$ is the $j$th row of $A$}
            
            $\mathcal{X}^{(j + n_x)} = x_k^i - \left(\sqrt{(n_x + \lambda)P_k^{x^i}}\right)_j$
       }
        
        $\mathcal{Y}^{(j)} = h(\mathcal{X}^{(j)}),$ \qquad \qquad $j = 0 ,1, \ldots  2n_x$ 
            
        $w_0 = \frac{\lambda}{n_x + \lambda}$, \qquad \qquad $w_j = \frac{1}{2(n_x + \lambda)}$  \qquad $j = 1, 2 \ldots  2n_x$

        $y_k^i = \sum_{j=0}^{2n_x} w_j^{(m)} \mathcal{Y}^{(j)}$  
        
        $P_k^{y^i} = \sum_{j=0}^{2n_x} w_j^{(c)} (\mathcal{Y}^{(j)} - y_k^i)(\mathcal{Y}^{(j)} - y_k^i)^T + R_k$ 
    
        $\tilde{w}_k^i = w_{k-1}^i \cdot p(y_k | x_k^i)$
   }
    $w_k^i = \frac{\tilde{w}_k^i}{\sum_{j=1}^{N} \tilde{w}_k^j}$

    $\hat{x}_k = \sum_{i=1}^{N} w_k^i x_k^i$\tcc*[r]{Estimation of the best solution (MMSE)}

    $\hat{y}_k = \sum_{i=1}^{N} w_k^i y_k^i$\tcc*[r]{Estimation of the best objective value (MMSE)}

    $P_k^{xx} = \sum_{i=1}^{N} w_k^i (x_k^i - \hat{x}_k)(x_k^i - \hat{x}_k)^T$

    $P_k^{yy} = \sum_{i=1}^{N} w_k^i (y_k^i - \hat{y}_k)(y_k^i - \hat{y}_k)^T$

    \If{$P_k^{xx} \prec P_{\text{min}}^{xx}$ \textbf{or} $P_k^{yy} \prec P_{\text{min}}^{yy}$}{
        \textbf{break} \tcc*[r]{Convergence check}
    }
}
\end{algorithm}

The Minimum Mean Square Error (MMSE) estimates for the state and objective value at iteration \( k \) are computed as,
\[
\hat{{x}}_k = \sum_{i=1}^N w_k^i {x}_k^i, \quad \hat{{y}}_k = \sum_{i=1}^N w_k^i {y}_k^i.
\]
In the context of the proposed PFO algorithm, the posterior distribution over the optimization variable is implicitly represented by a weighted ensemble of particles. As iterations progress, the accumulation of informative measurements induces posterior contraction, i.e., the posterior distribution increasingly concentrates around the minimum candidate.
This contraction is quantitatively captured by the empirical covariance matrix,
\[
P_k^{xx} = \sum_{i=1}^N w_k^i (x_k^i - \hat{x}_k)(x_k^i - \hat{x}_k)^\top.
\]
A similar covariance $P_k^{yy}$ is computed in the objective space, as,
\[
{P}_k^{{yy}} = \sum_{i=1}^N w_k^i ({y}_k^i - \hat{{y}}_k)({y}_k^i - \hat{{y}}_k)^\top.
\]
The decay of these covariances over time serves as an indicator of convergence. As $P_k^{xx} \to 0$, the particles collapse toward the estimated global minimizer, reflecting reduced uncertainty and enhanced confidence in the solution. This behavior aligns with the theory of posterior contraction in Bayesian inference, where the posterior mass increasingly concentrates in arbitrarily small neighborhoods of the true solution under mild assumptions on the noise model and proposal distribution \cite{ghosal2017fundamentals}.
Accordingly, if either of these covariances falls below a predefined bound,
the algorithm terminates, concluding that the particles have sufficiently concentrated around a global minimum candidate.

\begin{remark}
Resampling may play a critical role in standard particle filtering frameworks to mitigate weight degeneracy \cite{doucet2000sequential}. Over time, without resampling, the weight distribution tends to collapse, leaving only a few particles with significant weights. Incorporating resampling mechanisms such as systematic or residual resampling can enhance numerical stability and ensure better particle diversity. Accordingly, incorporation of the PFO algorithm with resampling can further improve robustness and convergence behavior.
\end{remark}

\subsection{Modeling the Transitional Prior Using Covariance Ellipsoids}\label{subsec:transitional}

In PFO, the transitional prior $p(x_k^i|x_{k-1}^i)$ can be any arbitrary density function that satisfies the following set of conditions:
\begin{enumerate}
    \item The transitional prior should encourage particles near the best estimated solution for exploitation.
    \item It should encourage particles far from the best estimated solution for exploration.
    \item It must incorporate sufficient stochasticity or noise to maintain a non-zero probability of reaching unexplored regions.
    \item The prior should ensure adequate coverage of the search space across iterations, maintaining smooth particle transitions.
\end{enumerate}
Since the UT is employed to propagate the true mean and covariance through a nonlinear transformation to the likelihood space, we propose constructing covariance ellipsoids based on the sigma points and use the geometry of the corresponding covariance matrix for directional updates of the optimization algorithm. Let the augmented data matrix be defined as,
\begin{align}
\xi = \left[\begin{array}{cc}
\mathcal{X}^0 & \mathcal{Y}^0\\
\mathcal{X}^1 & \mathcal{Y}^1\\
\vdots & \vdots\\
\mathcal{X}^{2n_x} & \mathcal{Y}^{2n_x}\\
\end{array}\right] \in \mathbb{R}^{(2n_x+1)\times (n+1)},
\end{align}
where $\mathcal{X}^j\in \mathds{R}^{1\times n}$ and $\mathcal{Y}^j\in \mathds{R}$ for $j=0, \ldots ,2n_x$ are sigma points and their transformed values after UT, respectively.
To quantify the local geometry of the search around each particle, we compute the empirical covariance $C_{xy}$ of the augmented matrix $\xi$ as,
\begin{align*}
&\bar{\xi} = \frac{1}{2n_x + 1} \sum_{j=0}^{2n_x} \xi^j,\\
&C_{xy} = \frac{1}{2n_x} \left(\xi - \mathbf{1}_{2n_x+1} \cdot \bar{\xi}^\top\right)^\top \left(\xi - \mathbf{1}_{2n_x+1} \cdot \bar{\xi}^\top\right),
\end{align*}
where $\xi^j = [\mathcal{X}^j, \mathcal{Y}^j]$ is the $j$th row of $\xi$, and $\bar{\xi}$ is the mean of the sigma points in the augmented space.
Then the eigenvectors of $C_{xy}$ provide the principal axes of the joint ellipsoid formed in the generalized space (i.e., the space of optimization variables augmented with the objective function values) and the eigenvalues correspond to the variance along those directions \cite{jolliffe2011principal}. \\
Next, define the augmented direction vector for each particle, i.e., $d_{k-1}^i$ as,
\[
d_{k-1}^i = \begin{bmatrix}
\hat{x}_{k-1} - x_{k-1}^i \\
\hat{y}_{k-1} - y_{k-1}^i
\end{bmatrix} \in \mathbb{R}^{n+1},
\]
where 
\(
(\hat{x}_{k-1}, \hat{y}_{k-1})
\)
represents the current weighted average of particles in the augmented space.
Define the step size as,
\begin{align}\label{equ:stepsize}
s^i_{k} = 
\begin{cases}
\|d_{k-1}^i\|, & \text{} d_{k-1}^{i^\top} C^{-1} {d_{k-1}^i} \leq 1 \\
\displaystyle {\dfrac{\|d_{k-1}^i\|}{\|d_{k-1}^i\|_{C^{-1}_{xy}}}=\left(\dfrac{d_{k-1}^{i^\top} d_{k-1}^i}{{d_{k-1}^i}^\top C_{xy}^{-1}d_{k-1}^i}\right)^{\frac{1}{2}}}. & \text{} d_{k-1}^{i^\top} C^{-1} {d_{k-1}^i} > 1
\end{cases}
\end{align}

To ensure a non-zero probability for unvisited places in the search space, a zero-mean noise with normal distribution and covariance $Q$ is added to the transitional prior. Therefore, the particle's local update function or transitional prior can be written as:
\begin{align}
x_k^i\sim N\left(x_{k-1}^i+s_{k}^i \left(\hat{x}_{k-1} - x_{k-1}^i\right) ,Q\right)
\end{align}

To investigate the effectiveness of the proposed transitional posterior, let,

\[
Z_i = [z^i_1, z^i_2, \ldots, z^i_{n+1}] \in \mathbb{R}^{(n_x+1) \times (n+1)},
\]
be the matrix of eigenvectors of the local covariance matrix $C_{xy}$, and ,
\[
\Lambda_i = \mathrm{diag}(\lambda^i_1, \lambda^i_2, \ldots, \lambda^i_{n+1}),
\]
be the diagonal matrix of the corresponding eigenvalues. Then $C_{xy}=Z_i \Lambda_i Z_i^\top$ and
${d}_{k-1}^{i,\mathrm{proj}} = Z_i^T {d}_{k-1}^i$ gives the projection of the direction ${d}_{k-1}^i$ on the eigenvector basis of the covariance ellipsoid. According to the proposed update rule, if $\begin{bmatrix}
\hat{x}_{k-1}  \\
\hat{y}_{k-1}
\end{bmatrix}$ is not within the ellipsoid, then the step-size is,
\begin{align*}
    s^i_k &=\left(\dfrac{d_{k-1}^{i^\top} d_{k-1}^i}{d_{k-1}^{i^\top} C_{xy}^{-1}d_{k-1}^i}\right)^{\frac{1}{2}} = \left(\dfrac{d_{k-1}^{i^\top} d_{k-1}^i}{d_{k-1}^{i^\top} Z_i \Lambda_i^{-1} Z_i^\top d_{k-1}^i}\right)^{\frac{1}{2}}\\
    &  = 
    \left(\dfrac{d_{k-1}^{i^\top} d_{k-1}^i}{{d}_{k-1}^{{i,\mathrm{proj}}^\top} \Lambda_i^{-1} {d}_{k-1}^{i,\mathrm{proj}}}\right)^{\frac{1}{2}}=
    \left(\dfrac{d_{k-1}^{i^\top} d_{k-1}^i}{\sum_{j=1}^{n+1} \frac{({d}_{k-1}^{i,\mathrm{proj}}(j))^2}{\lambda^i_j}}\right)^{\frac{1}{2}}.
\end{align*}
This interpretation reveals that directions aligned with low-variance eigenvectors are penalized more heavily, effectively regularizing against stepping too far in directions of uncertainty.

As a special case, when the direction vector \( d_{k-1}^i \) lies outside the ellipsoid and aligns with one of the eigenvectors of \( C_{xy} \), say \( z_j \), we can write $d_{k-1}^i = \alpha z_j,$
for some scalar \( \alpha \). Then $s_k^i$ simplifies to,

\begin{equation*}
    s_k^i = \left( \dfrac{\|d_{k-1}^i\|^2}{\|d_{k-1}^i\|_{C^{-1}_{xy}}^2} \right)^{1/2}
    = \left( \dfrac{\alpha^2}{\alpha^2 \cdot \frac{1}{\lambda^i_j}} \right)^{1/2}
    = \sqrt{\lambda^i_j}
\end{equation*}
This behavior shows how the algorithm makes cautious updates along sensitive directions and bolder moves along directions with higher local uncertainty. In practical optimization algorithms, a step size gain $\gamma$ is typically applied to $s_k^i$. This introduces a degree of freedom, allowing for control over the explore-exploit balance throughout the optimization process. \\
Moreover, this covariance-dependent update rule naturally satisfies both the exploitation (condition (i)) and exploration (condition (ii)) requirements of the transitional prior: particles located closer to the current estimated minimum exhibit tighter covariance ellipsoids (i.e., smaller eigenvalues), thereby promoting local exploitation. In contrast, particles farther from the minimum tend to have broader ellipsoids (i.e., larger eigenvalues), encouraging more global exploration. Moreover, zero-mean noise with normal distribution and covariance $Q$ in the transitional prior ensures satisfaction of condition (iii).
Finally, by carefully modulating the spread and direction of these ellipsoids during the resampling or prediction step, the proposed method can span the search space smoothly across iterations, thereby addressing Condition (iv) concerning continuous and comprehensive space coverage.\\

\begin{figure}[H]
      \centering
      \includegraphics[width=.4\columnwidth]{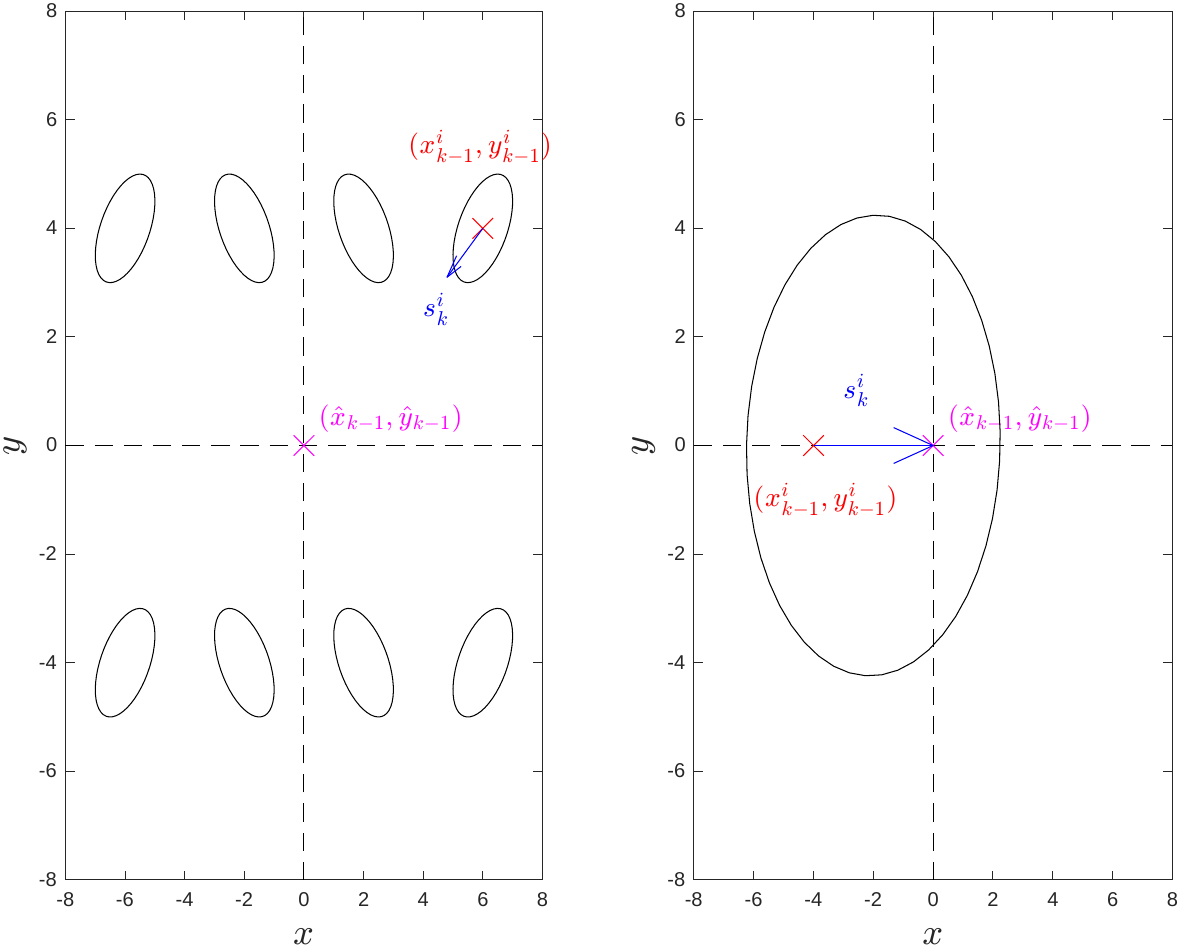}\\
     \caption{Local update function or state transitional prior function}\label{fig:localUpdate}
  \end{figure}

Figure \ref{fig:localUpdate} displays the proposed method for particle motions in 2D (it's readily generalizable for higher dimensions). As this figure demonstrates, based on the position of the ellipsoid mean value (center) with respect to the estimated solution, the step size is chosen as by $s_k^i$ for particle-$i$. Then, whether the ellipsoid contains the estimated solution or not, the step size is adjusted according to \eqref{equ:stepsize}. In both scenarios, the direction of particle motion will be toward the estimated solution.

\begin{figure}[H]
      \centering
      \subfloat[]{\label{}\includegraphics[width=0.32\columnwidth]{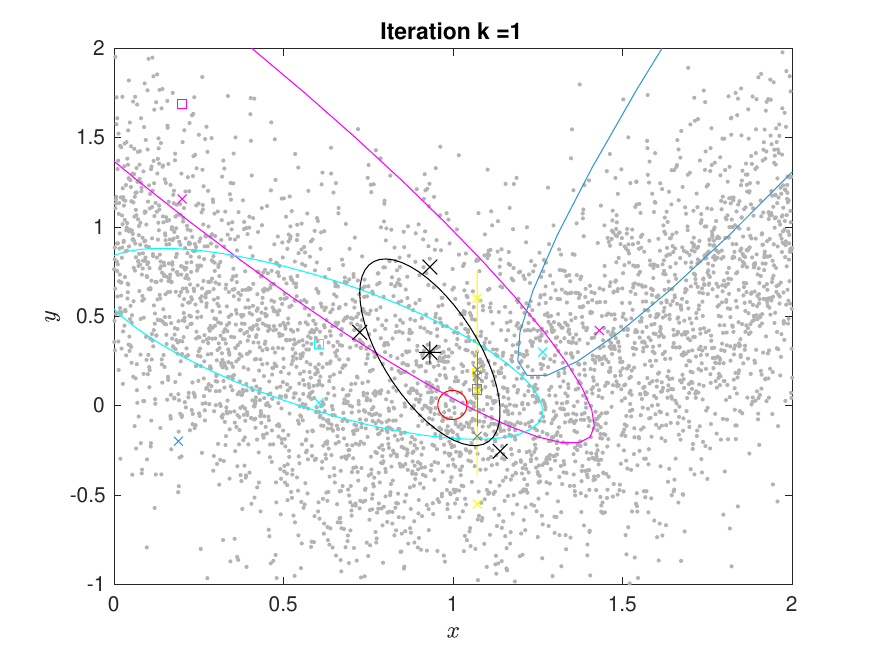}}
      \subfloat[]{\label{}\includegraphics[width=0.32\columnwidth]{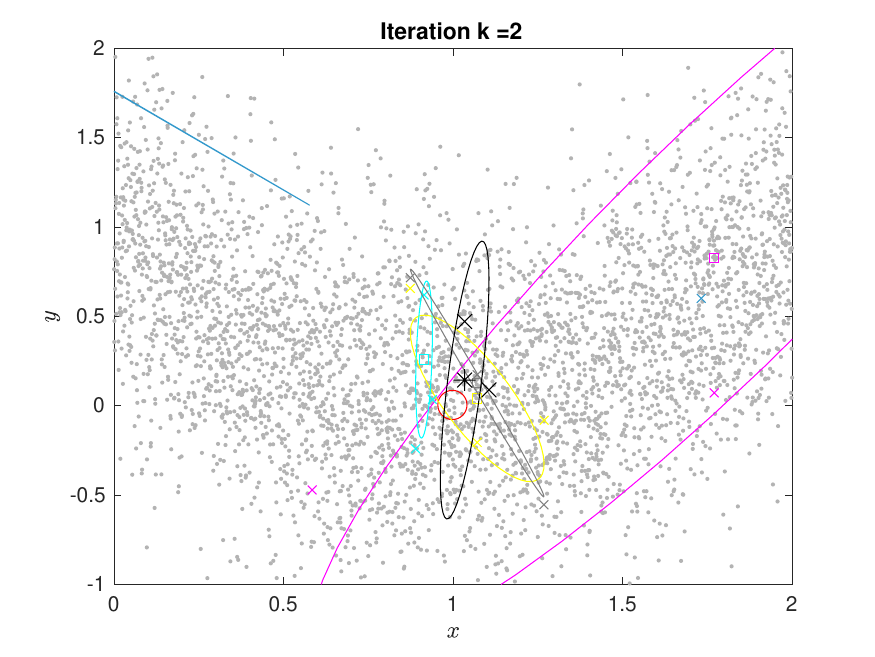}}
      \subfloat[]{\label{}\includegraphics[width=0.32\columnwidth]{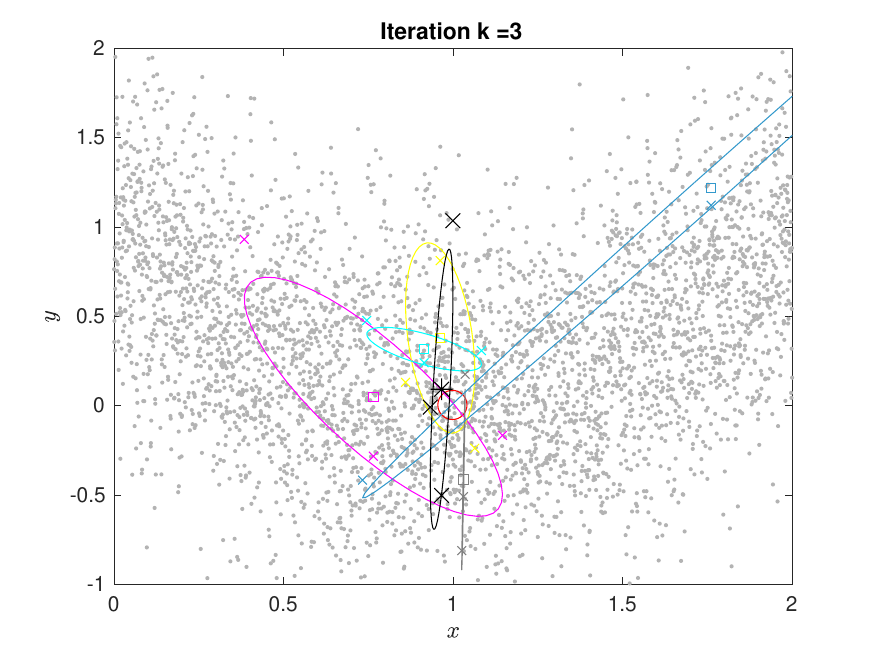}}\\
      \subfloat[]{\label{}\includegraphics[width=0.32\columnwidth]{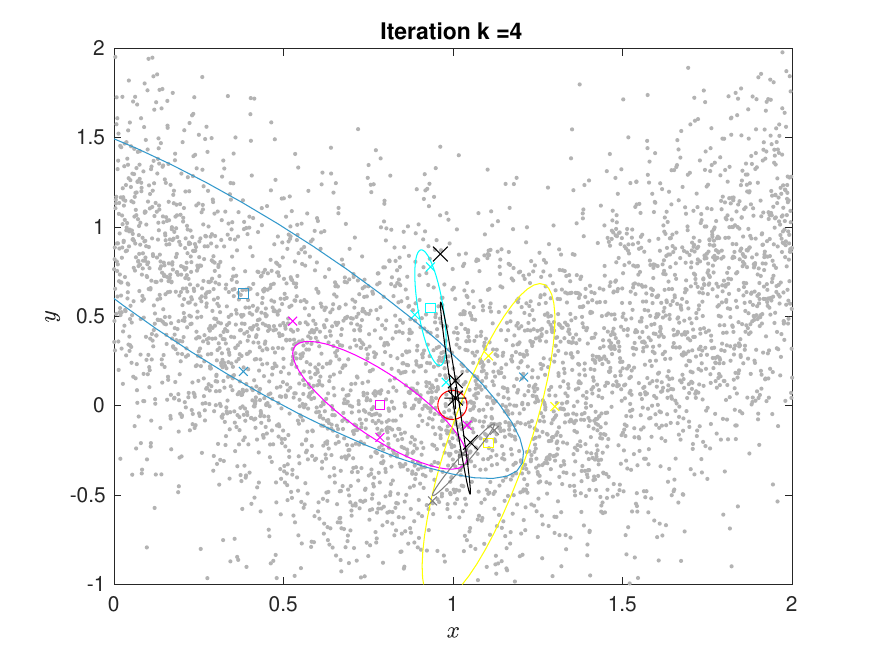}}
      \subfloat[]{\label{}\includegraphics[width=0.32\columnwidth]{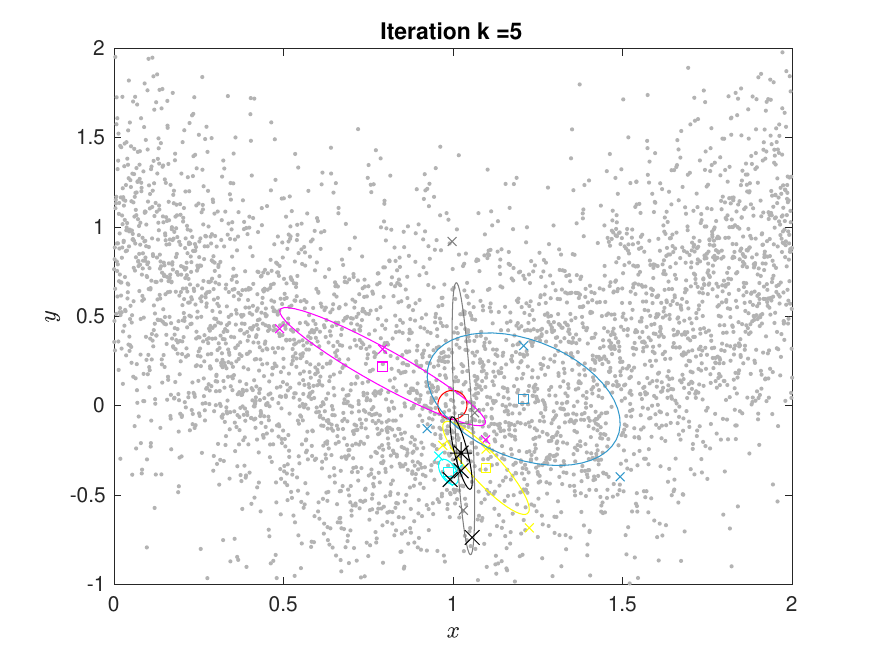}}
      \subfloat[]{\label{}\includegraphics[width=0.32\columnwidth]{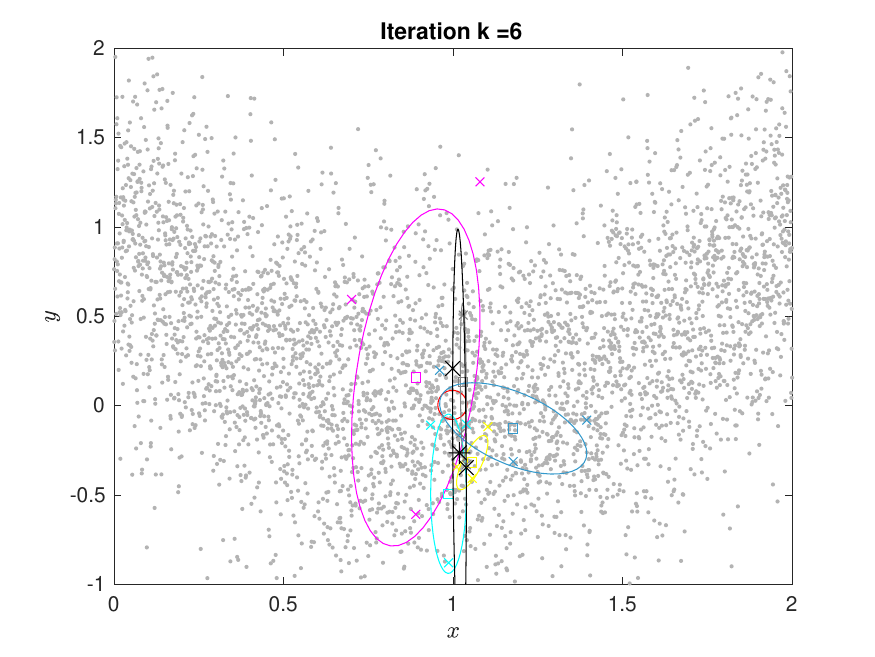}}\\
      \subfloat[]{\label{}\includegraphics[width=0.32\columnwidth]{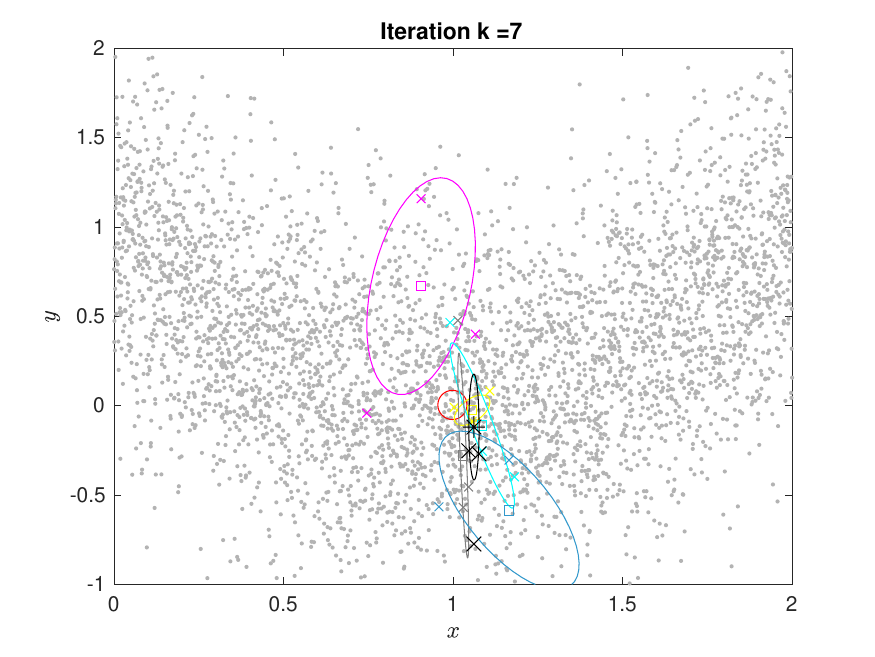}}
      \subfloat[]{\label{}\includegraphics[width=0.32\columnwidth]{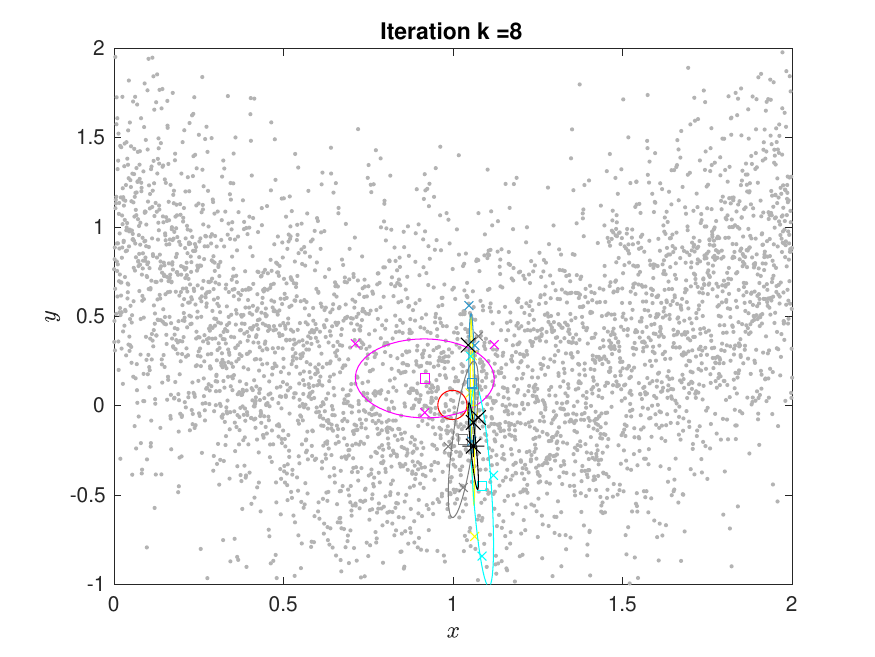}}
      \subfloat[]{\label{}\includegraphics[width=0.32\columnwidth]{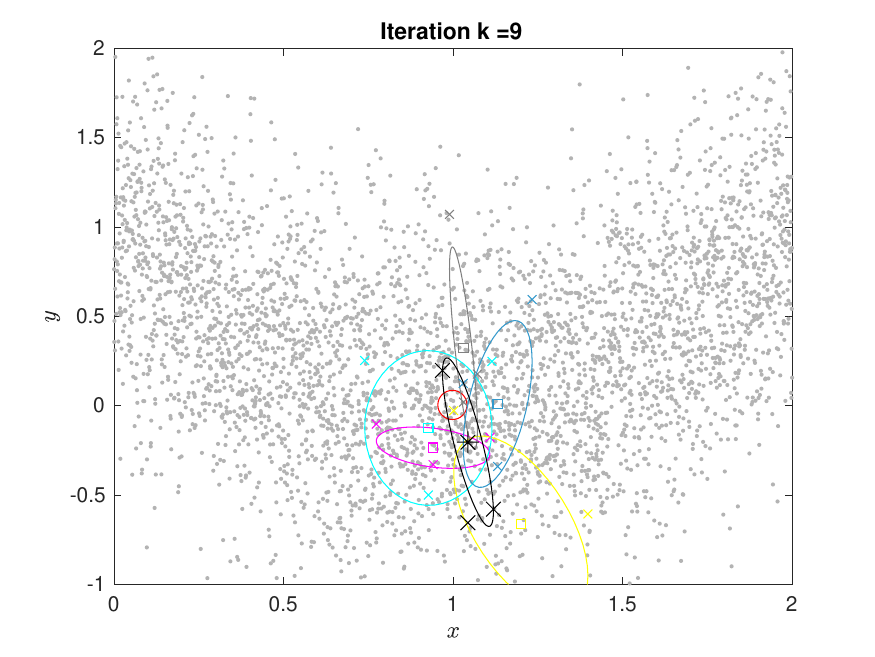}}
     \caption{Particle motions in iterations}\label{fig:motions}
  \end{figure}

To provide empirical insight into the particle motions under the command of this transitional prior, Figure \ref{fig:motions} illustrates particle covariance ellipsoids, sigma points ($\times$ markers), and their mean value ($\square$ markers) for 9 steps in an example problem for $N=5$ (problem number 2 in Table \ref{tab:example_func}). The $\ast$ marker shows the best solution found at iterations with the corresponding sigma point and covariance ellipsoid in black. The red circle represents the actual minima.

Roughly speaking, these plots show that the magnet and blue particles are exploring the area until iteration 7 and then settle down near the estimated minima. At the same time, the cyan, gray, and yellow particles are exploiting to find a better solution near the actual best.

\subsection{Convergence and Confidence Level Analysis}

This section provides a rigorous analysis of the convergence properties and associated confidence measures of the PFO algorithm. PFO operates within a Bayesian inference framework, combining importance sampling with mechanisms akin to posterior contraction to perform global stochastic optimization under noisy conditions.

\subsubsection{Posterior Contraction}

Let $x^\ast \in \mathbb{R}^n$ denote the global minimizer of the objective function $h(x)$, and let $p(x_k \mid y_{1:k})$ represent the posterior distribution over the candidate solution at iteration $k$, conditioned on all noisy measurements $y_{1:k}$. As the number of iterations increases and more informative measurements are incorporated, the posterior contracts around $x^\ast$. Formally, for any $\epsilon_x > 0$,
\begin{align}
    \mathbb{P}\left( \| x_k - x^\ast \| > \epsilon_x \mid y_{1:k} \right) \to 0, \;\; \text{as } k \to \infty,
\end{align}
under mild assumptions on the noise model and the informativeness of the transition prior (as discussed in sub-section~\ref{subsec:transitional}). This is consistent with general results on posterior consistency and contraction in nonparametric Bayesian inference~\cite{ghosal2017fundamentals}.

The rate of contraction can be empirically monitored via the covariance matrices $P_k^{xx}$ and $P_k^{yy}$, which capture the spread of particles in the input and objective spaces. As $P_k^{xx} \to 0$, particles collapse around the estimated minimizer, reflecting reduced uncertainty and increased confidence.

\subsubsection{Importance Sampling Approximation}

At each iteration, the posterior is approximated via importance sampling as:
\begin{align}
    \hat{p}_N(x_k \mid y_{1:k}) = \sum_{i=1}^N w_k^i \delta_{x_k^i}(x_k),
\end{align}
where $\delta_{x_k^i}$ is the Dirac measure centered at particle $x_k^i$. Under appropriate regularity conditions (e.g., sufficient $N$, well-specified transition model, non-degenerate weights), this empirical posterior converges weakly to $p(x_k \mid y_{1:k})$ as $N \to \infty$~\cite{vo2005sequential}. Consequently, $\hat{x}_k$, the MMSE estimate, becomes a consistent estimator of the posterior mean, satisfying $\hat{x}_k \to x^\ast$ in probability. Resampling techniques such as systematic or residual resampling are adopted to mitigate weight degeneracy, improve approximation accuracy, and preserve diversity over iterations.

\subsubsection{Confidence Quantification and Stopping Criterion}

The matrix norms $\|P_k^{xx}\|$ and $\|P_k^{yy}\|$ are direct indicators of uncertainty. When the posterior distribution is approximately Gaussian in the neighborhood of $x^\ast$, one can define an elliptical $(1-\alpha)$-credible region as:
\begin{align}
    \mathcal{E}_x(\alpha) = \left\{ x \in \mathbb{R}^n : (x - \hat{x}_k)^\top (P_k^{xx})^{-1} (x - \hat{x}_k) \leq \chi^2_n(1 - \alpha) \right\},
\end{align}
where $\chi^2_n(1 - \alpha)$ is the quantile of the chi-squared distribution with $n$ degrees of freedom. Similar ellipsoids can be defined for the objective values via $P_k^{yy}$. Assuming diagonalize $P^{xx}_k\approx \sigma_x^2 I$, with the confidence level $1-\alpha\approx 1-\sigma_x$, the minimizer lies within radius $\sigma_x\sqrt{\chi_n^2(1-\sigma_x)}$ of the estimated solution $\hat{x}_k$.

The PFO algorithm terminates when both $P_k^{xx}$ and $P_k^{yy}$ fall below predefined thresholds $\sigma_x, \sigma_y > 0$, i.e.,
\begin{align}
    \|P_k^{xx} - P^{xx}_{\min}\| < \sigma_x, \quad \|P_k^{yy} - P^{yy}_{\min}\| < \sigma_y.
\end{align}
This ensures that the final estimates satisfy:
\begin{align}
    \|\hat{x}_k - x^\ast\| < \epsilon_x, \quad \|\hat{y}_k - y^\ast\| < \epsilon_y,
\end{align}
with confidence levels at least $(1 - \sigma_x)$ and $(1 - \sigma_y)$, respectively. This provides a theoretically grounded stopping criterion aligned with Bayesian credible regions.

\subsubsection{Joint Estimation and Bias Consideration}
Although the theoretical formulation, based on posterior contraction and importance sampling, rigorously specifies the expected convergence behavior, practical limitations inherent in nonlinear estimation methods such as PF and the UT introduce approximation errors and fundamental bias. As a result, one cannot arbitrarily select thresholds $\epsilon_x$ and $\epsilon_y$ independent of the associated confidence levels. Because, PFO algorithm performs joint estimation over both the input and objective values $(x, y)$ through recursive Bayesian updates. The UT is employed to propagate particles through the nonlinear objective $h(.)$, allowing for an estimation of the first two moments $(\hat{y}_k,P_k^{yy})$ of the output posterior. Thus, the algorithm performs a joint filtering task over $(x,y)$, and the convergence of $(\hat{x}_k,\hat{y}_k)$ to $(x^\ast,y^\ast)$ is governed by the contractive behavior of both posterior marginals under the designed dynamics.

In this context, the transition prior plays a dual role: It enables exploration-exploitation balance by controlling how particles move in the augmented space $(x,y)$, leveraging the local geometry of covariance ellipsoids around each particle; and, It ensures that particles which are closer to the estimated minimizer (low variance, high confidence) are refined locally, while particles farther away (higher variance, lower confidence) explore the objective landscape more broadly.

At iteration $k$, the empirical joint posterior is represented by:
\begin{align}
    \hat{p}_N(x_k, y_k \mid y_{1:k}) = \sum_{i=1}^N w_k^i \delta_{(x_k^i, y_k^i)}(x_k, y_k).
\end{align}
The MMSE estimates $(\hat{x}_k, \hat{y}_k)$, under standard assumptions (e.g., Lipschitz continuity of $h(.)$, bounded noise, non-degenerate weights) converges weakly in probability to $(x^\ast, y^\ast)$ as $N \to \infty$, provided the transition prior maintains sufficient coverage and the observation model is consistent. However, in practice, bias may emerge due to: (a) The UT underestimating higher-order moments in nonlinear regions; (b) Degeneracy in particle weights, especially under high noise or insufficient resampling; (c) Limited exploration in poorly tuned transition priors. Thus, it is possible for
\[
\lim_{k \to \infty} \|P_k^{xx}\| \to 0 \quad \text{does not necessarily imply} \quad \hat{x}_k \to x^\ast,
\]
unless the MMSE estimator is asymptotically unbiased. Mitigating this bias requires careful tuning of algorithm parameters and approximations such as resampling strategies, noise covariances, and adaptive UT scaling. These adjustments can be performed using a hyperparameter optimization algorithm to initially tune the parameters for specific problem instances, as demonstrated in the next section.

\section{Empirical Evaluation of PFO Performance and Robustness}\label{sec:eval}
\subsection{Example Set 1}
\begin{figure}[ht]
      \centering
      \subfloat[]{\label{}\includegraphics[width=0.45\columnwidth]{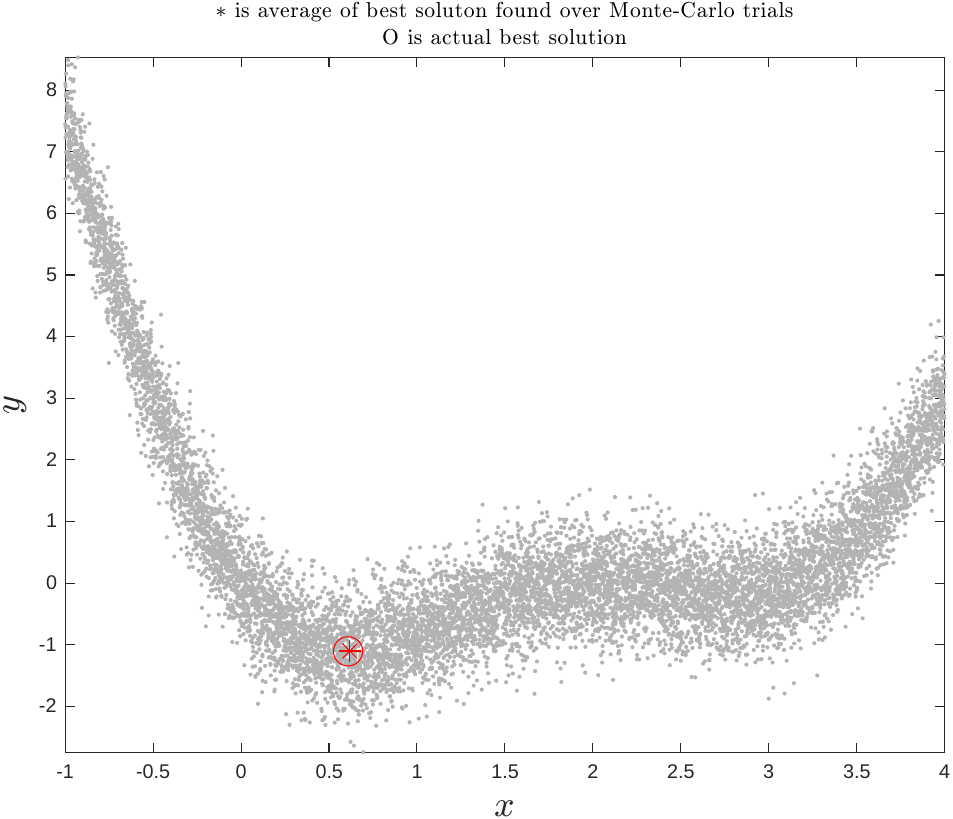}}
      \subfloat[]{\label{}\includegraphics[width=0.45\columnwidth]{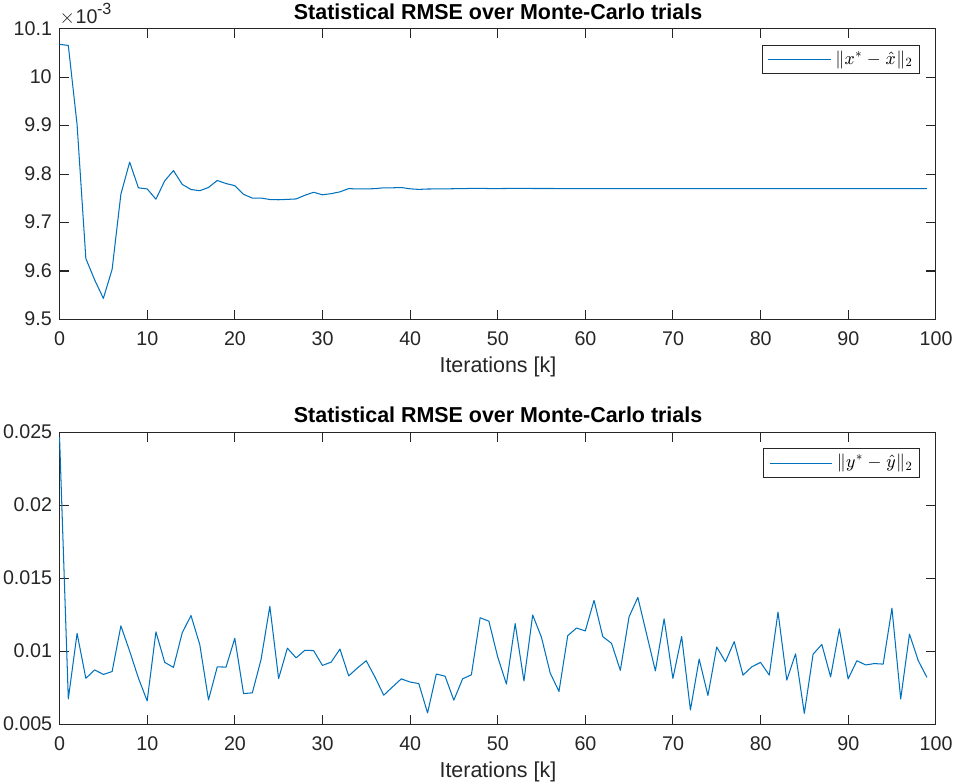}}\\
     \caption{Result of PFO for function $H_1(x)$}\label{fig:h1}
  \end{figure}

\begin{figure}[ht]
  \centering
  \subfloat[]{\label{}\includegraphics[width=0.45\columnwidth]{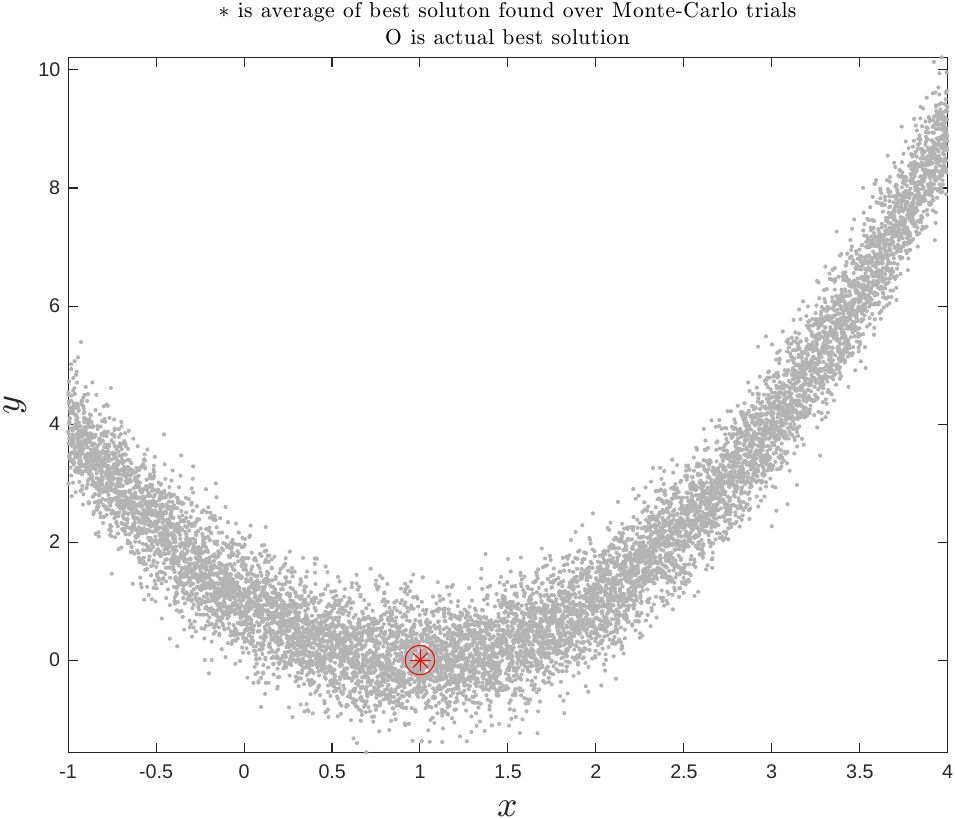}}
  \subfloat[]{\label{}\includegraphics[width=0.45\columnwidth]{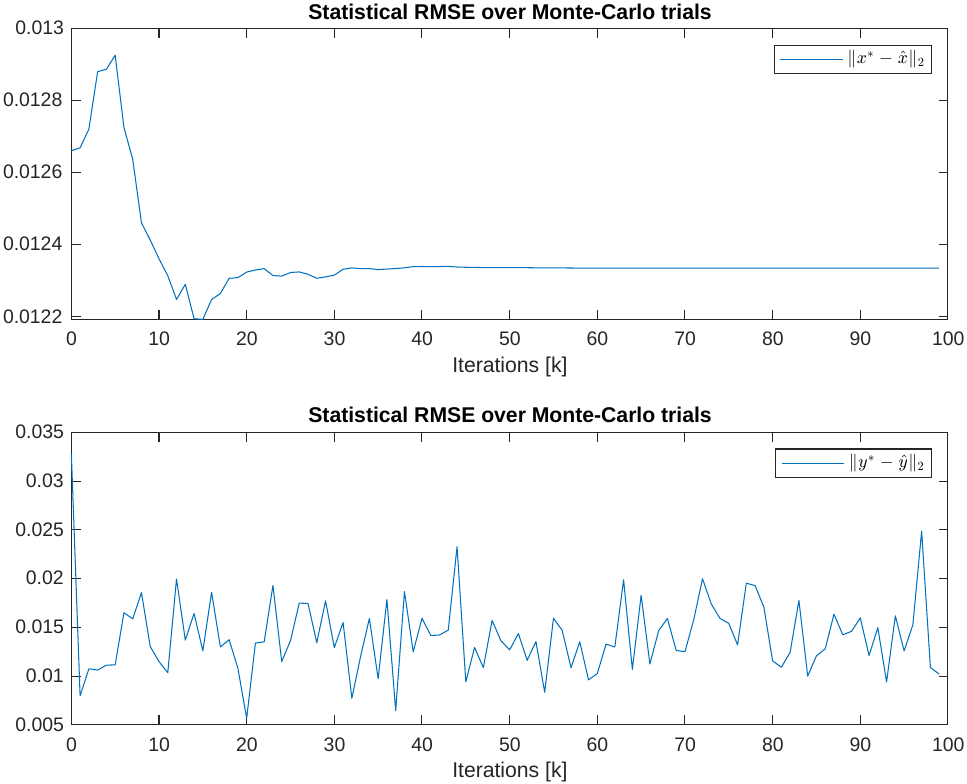}}\\
  \caption{Result of PFO for function $H_2(x)$}\label{fig:h2}
\end{figure}

\begin{figure}[ht]
  \centering
  \subfloat[]{\label{}\includegraphics[width=0.45\columnwidth]{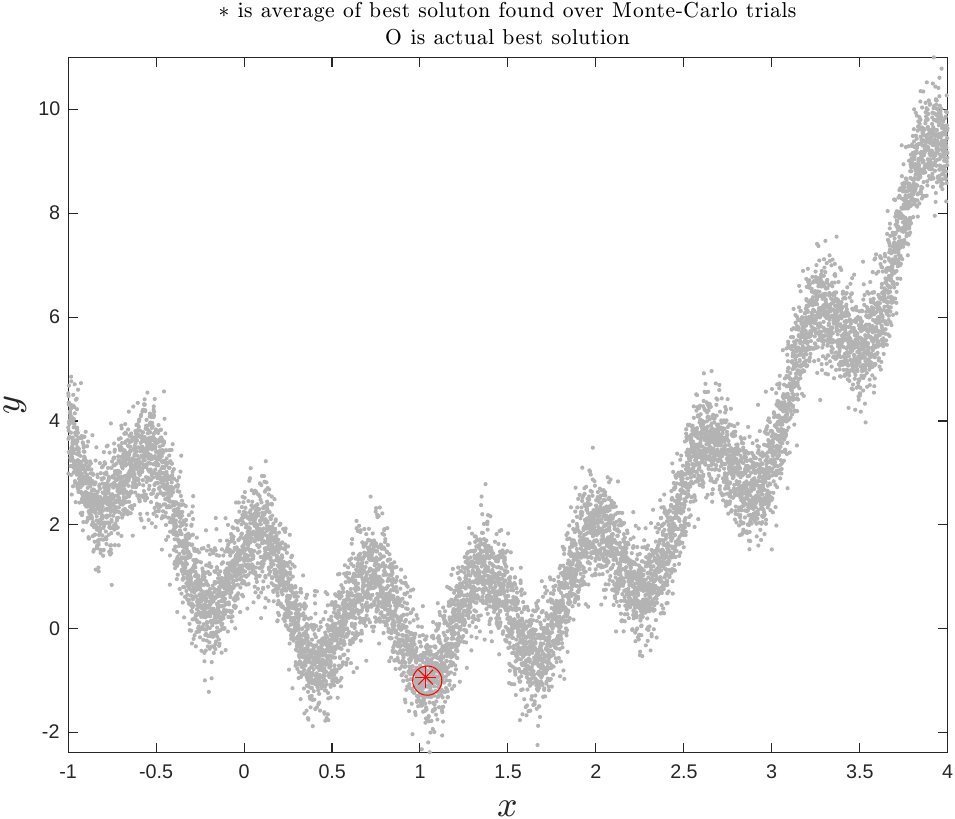}}
  \subfloat[]{\label{}\includegraphics[width=0.45\columnwidth]{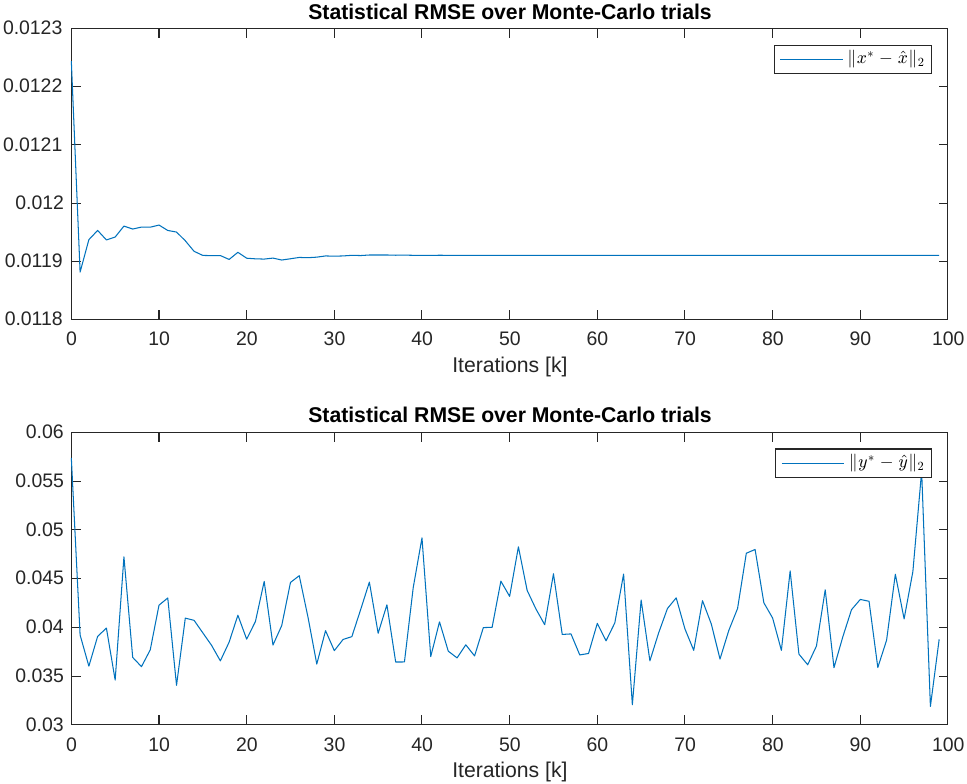}}\\
 \caption{Result of PFO for function $H_3(x)$}\label{fig:h3}
\end{figure}

\begin{figure}[ht]
  \centering
  \subfloat[]{\label{}\includegraphics[width=0.45\columnwidth]{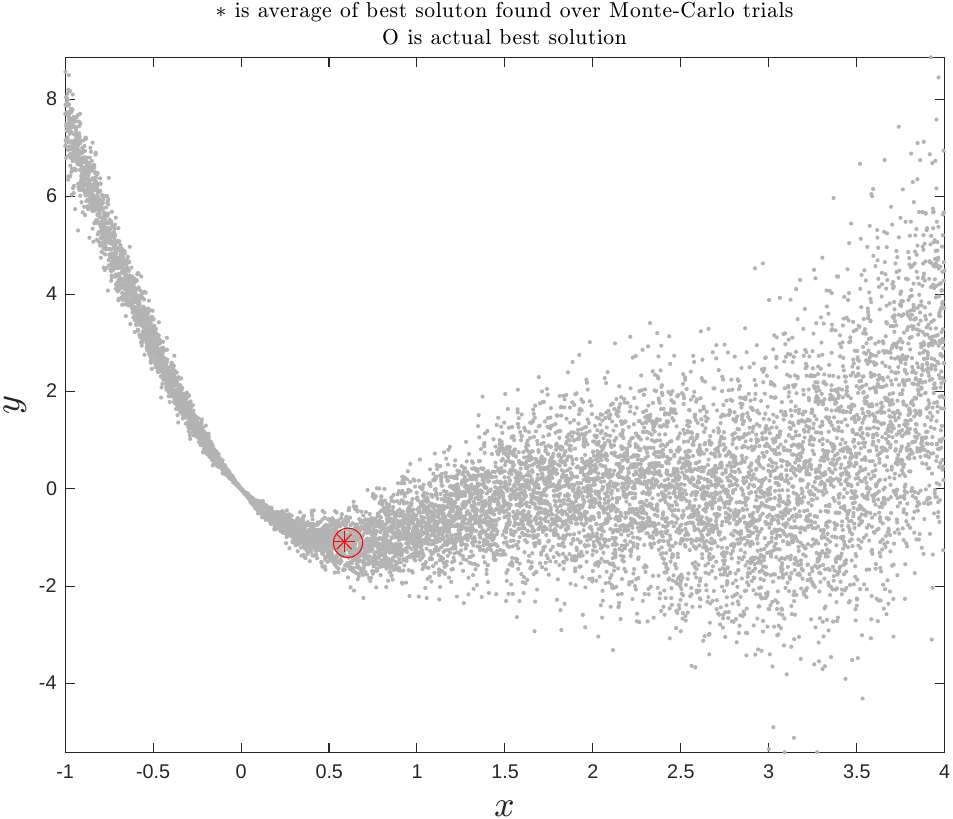}}
  \subfloat[]{\label{}\includegraphics[width=0.45\columnwidth]{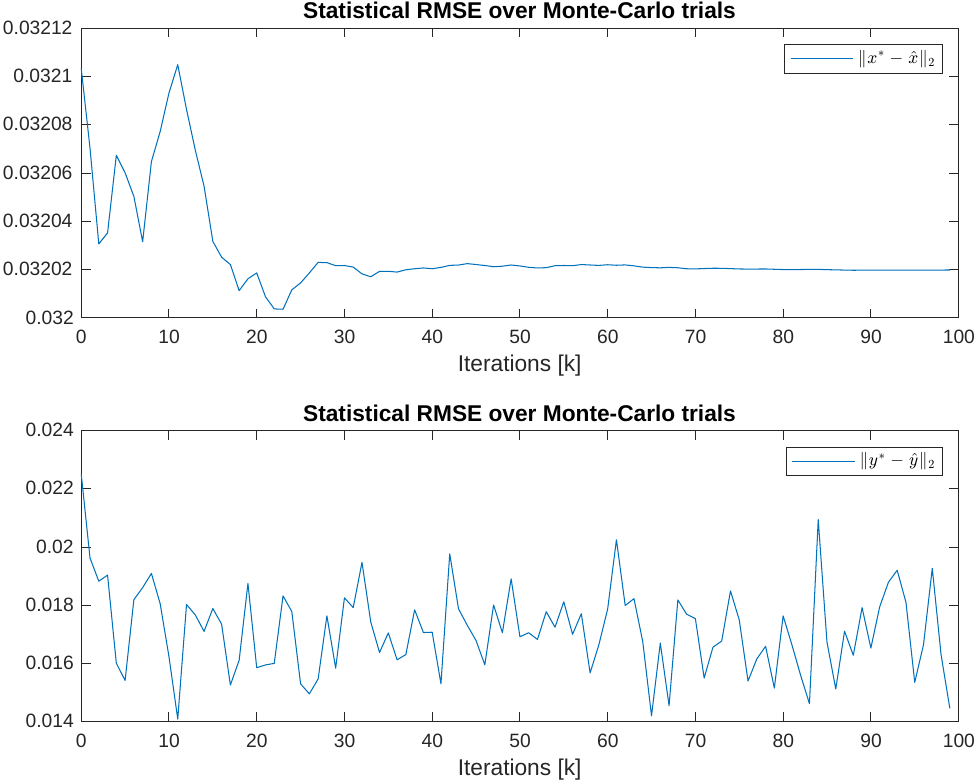}}\\
 \caption{Result of PFO for function $H_4(x)$}\label{fig:h4}
\end{figure}

Performance and robustness of the proposed PFO algorithm are tested using the example functions presented in Table \ref{tab:example_func}. The robustness is checked through 10 Monte Carlo trials. The algorithm parameters presented in this table are tuned using random sampling over predefined set. Figures \ref{fig:h1} to \ref{fig:h4} display the problem's data spread and the best-found solution in (a), while (b) demonstrates the statistical RMSE for each step over the Monte Carlo trials. The plots indicate high confidence in finding the global minima within pre-determined uncertainty bounds. In all optimizations the step gain is identity, i.e. $\gamma=1$.

\begin{table*}
  \centering
  \begin{tabular}{lllllllll}
    \hline
    Functions & $k_{max}$ & $N$ & $N_{thr}$ &  $\lambda$ & $Q$ & $\sigma_x$ & $\sigma_y$ & $R$\\\hline\hline
    $H_1(x) = -\sin(x)(x-2)^2+v$ & 100 & 100 & $N/2$ & 1 & 1e-8 & 1e-5 & 1e-5 & 0.5\\
    $H_2(x) = (x-1)^2+v$ & 100 & 50 &  $N/2$ & 1&  1e-8 & 1e-5 & 1e-5 & 0.5 \\
    $H_3(x) = (x-1)^2+\cos(10(x-0.1))+v$ & 100 & 50 & $N/2$ & 1 & 1e-8 & 1e-5 & 1e-5 & 0.5 \\
    $H_4(x) = -\sin(x)(x-2)^2+vx$ & 100 & 100 & $N/2$ & 2 & 1e-7 & 1e-5 & 1e-5 & 0.5 \\
    \hline
  \end{tabular}
  \caption{Example function set's parameter table}\label{tab:example_func}
\end{table*}


\subsection{Example Set 2: CEC 2005 Benchmark}
In this study, the PFO algorithm is tested using benchmark functions 1 and 4 from the 2005 CEC technical report~\cite{suganthan2005problem}, denoted as \( f_1(x) \) and \( f_4(x) \), respectively. This report includes 25 benchmark functions and evaluates the performance of several real-parameter optimization algorithms. The benchmark function is widely accepted as a baseline for newer function sets introduced in CEC 2013~\cite{li2013benchmark} and CEC 2022~\cite{luo2022benchmark}, and the PFO has been used for comparison with many functions in those sets. In this work, only simple functions are selected and reported, as the focus is on evaluating the qualitative behavior of the PFO algorithm rather than conducting an extensive benchmarking study.

Although these functions do not fall within the class of noisy target functions considered in this paper, artificial Gaussian noise can be added to their outputs to simulate noisy conditions. Alternatively, the PFO algorithm can be compared with other methods under zero-noise conditions, with some minor modifications. However, it is important to note that due to the computational overhead associated with handling uncertainty, the PFO algorithm is not ideally suited for noise-free problems and may exhibit longer runtimes compared to other heuristic methods such as PSO.

\begin{figure}[ht]
      \centering
      \subfloat[]{\label{}\includegraphics[width=0.45\columnwidth]{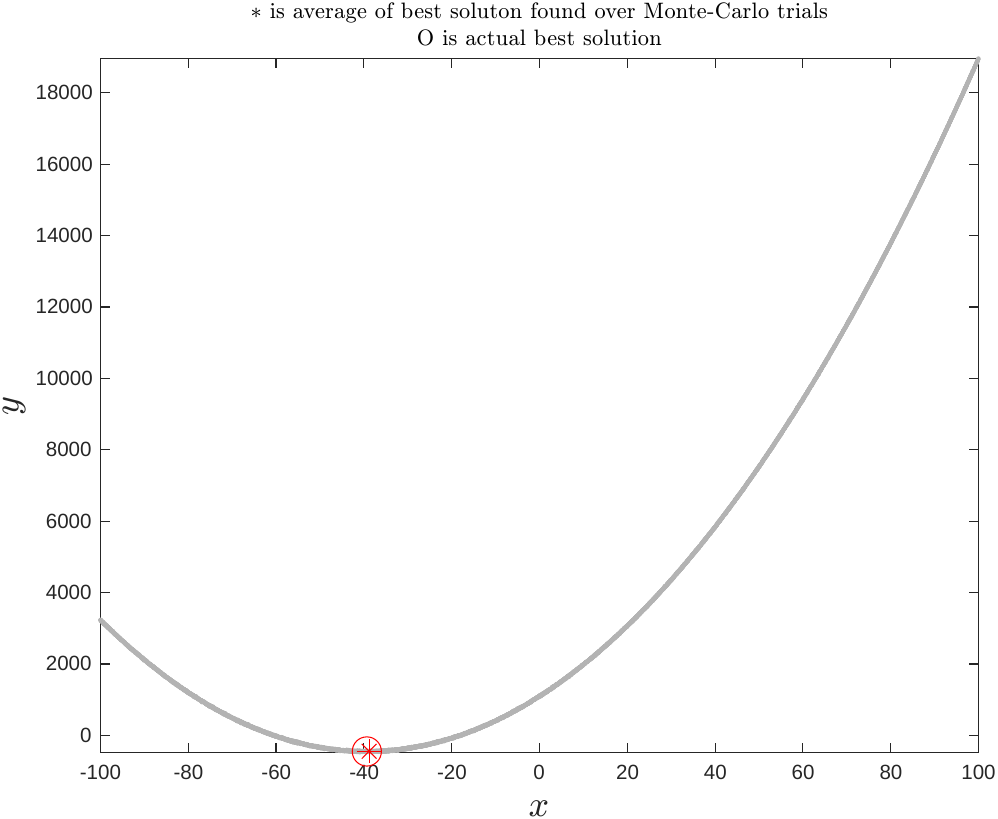}}
      \subfloat[]{\label{}\includegraphics[width=0.45\columnwidth]{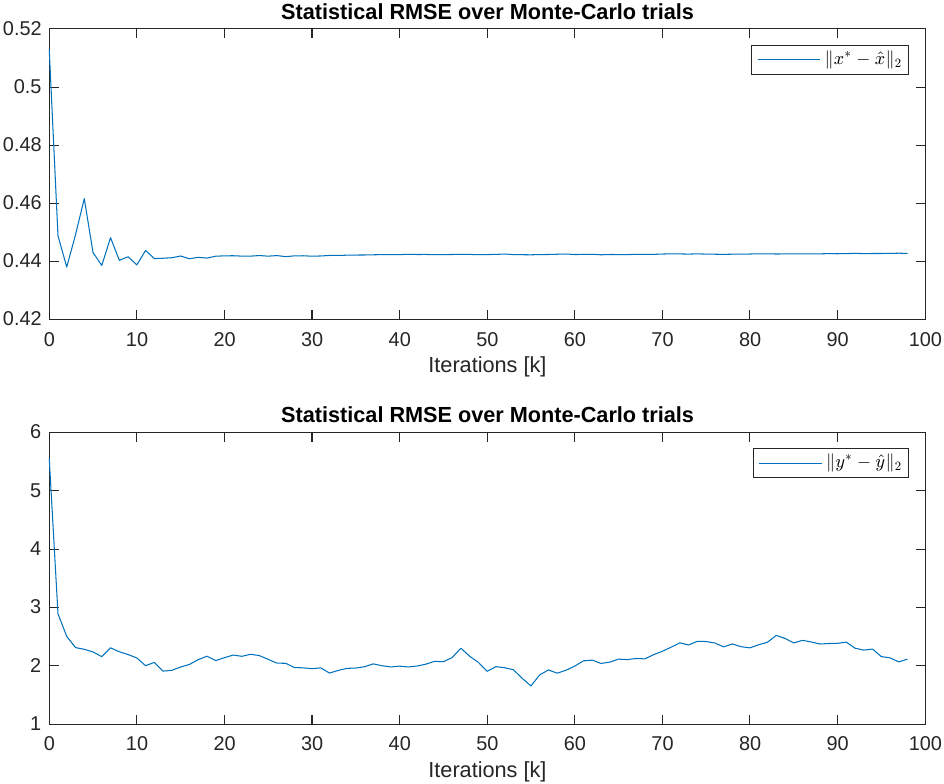}}\\
     \caption{Result of PFO for function $H_5(x)$}\label{fig:f1}
  \end{figure}

\begin{figure}[ht]
  \centering
  \subfloat[]{\label{}\includegraphics[width=0.45\columnwidth]{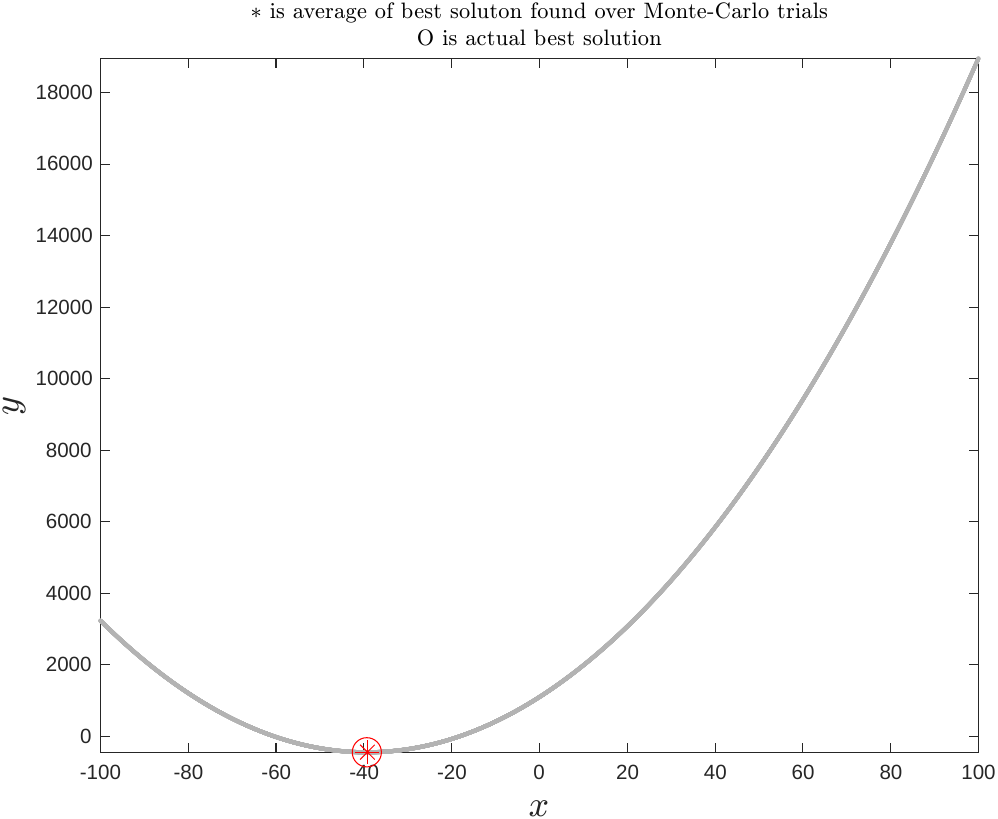}}
  \subfloat[]{\label{}\includegraphics[width=0.45\columnwidth]{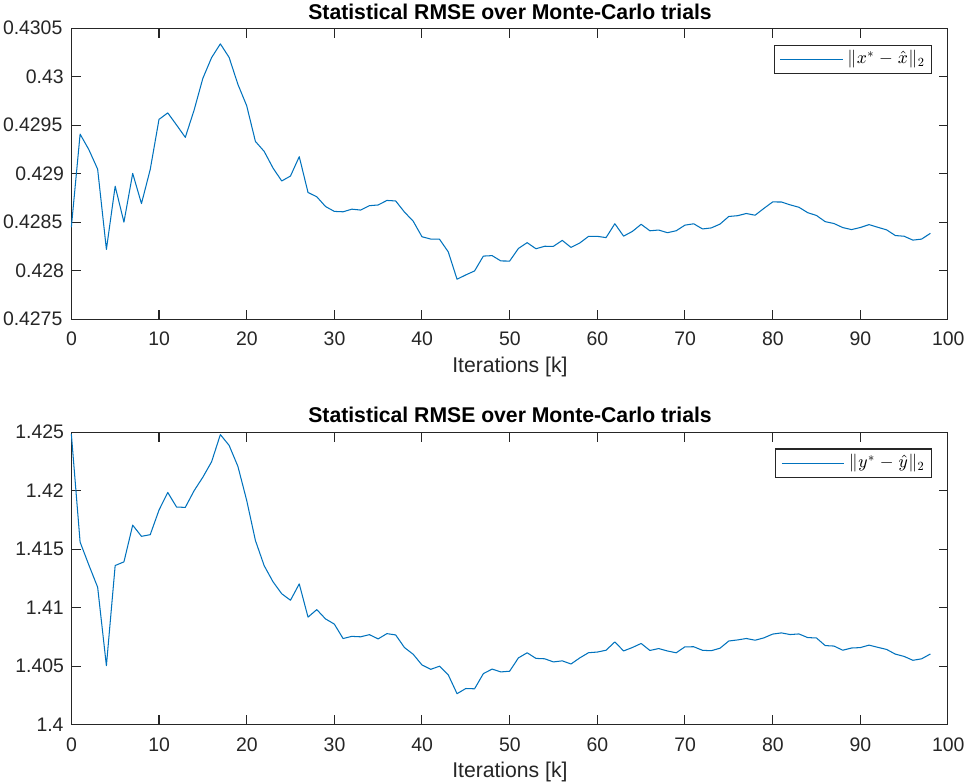}}\\
  \caption{Result of PFO for function $H_6(x)$}\label{fig:f2}
\end{figure}

\begin{figure}[ht]
  \centering
  \subfloat[]{\label{}\includegraphics[width=0.45\columnwidth]{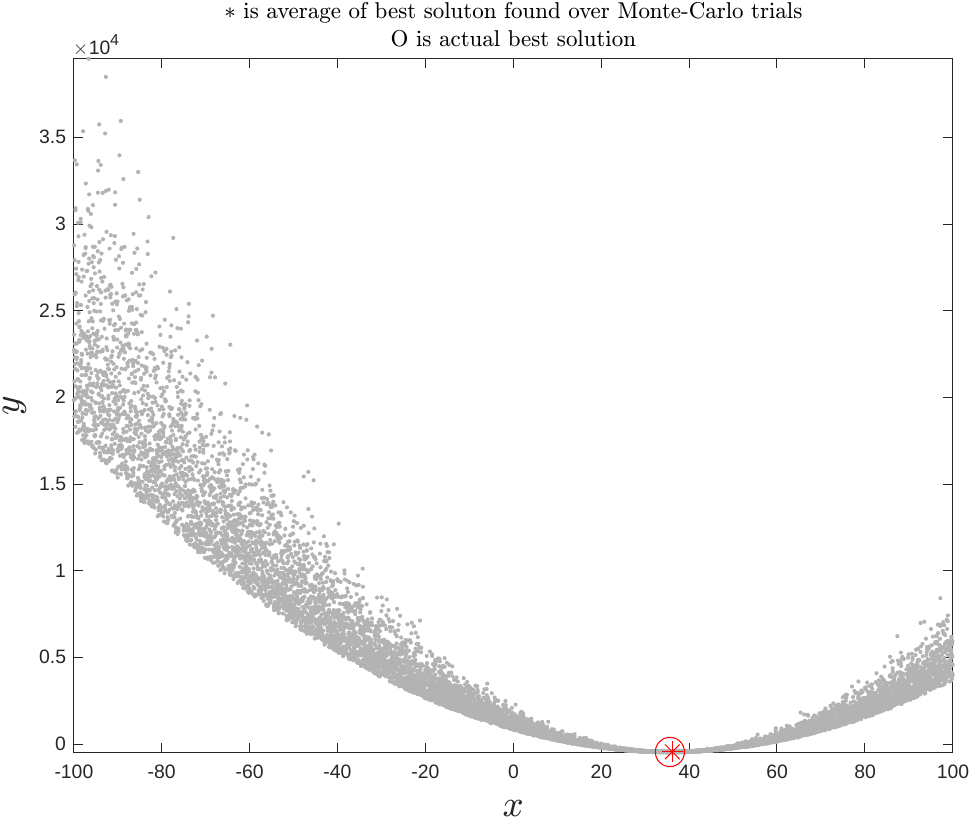}}
  \subfloat[]{\label{}\includegraphics[width=0.45\columnwidth]{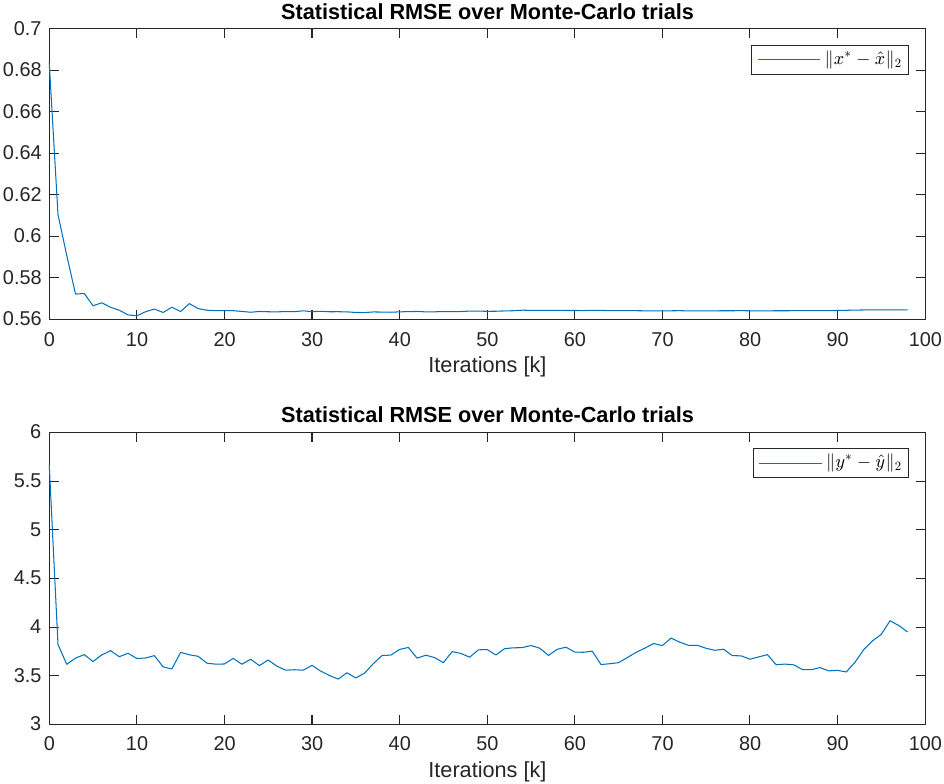}}\\
 \caption{Result of PFO for function $H_7(x)$}\label{fig:f3}
\end{figure}

\begin{figure}[ht]
  \centering
  \subfloat[]{\label{}\includegraphics[width=0.45\columnwidth]{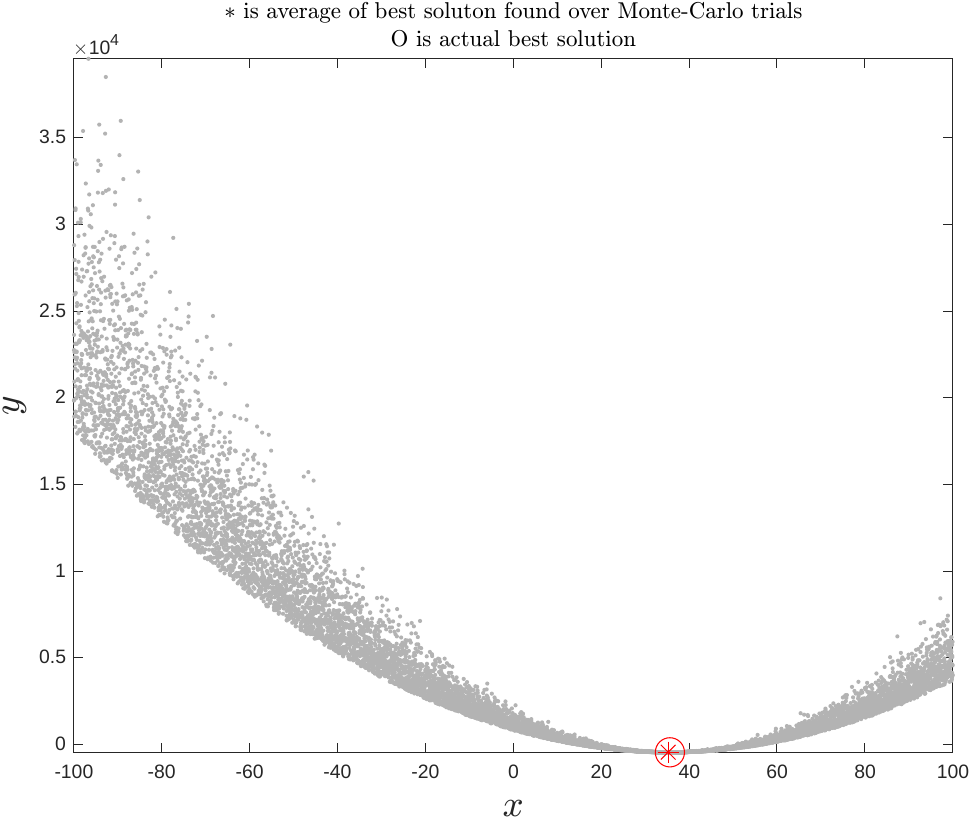}}
  \subfloat[]{\label{}\includegraphics[width=0.45\columnwidth]{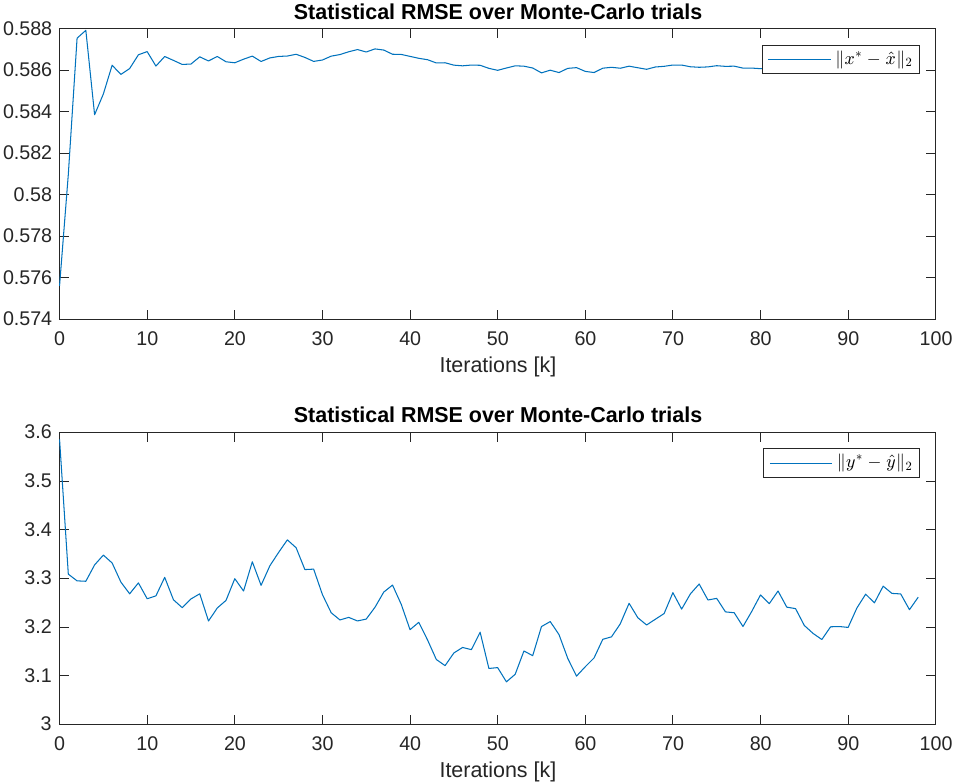}}\\
 \caption{Result of PFO for function $H_8(x)$}\label{fig:f4}
\end{figure}

For this reason, the comparison is restricted to a low-dimensional setting, specifically \( D = 1 \). The parameter settings used in this comparison are provided in Table~\ref{tab:example_func_cec}. Figures~\ref{fig:f1} to \ref{fig:f4} illustrate the performance of the PFO algorithm. Robustness is assessed through 25 Monte Carlo trials.  The algorithm parameters presented in this table are tuned using random sampling over predefined set.

\begin{table*}[ht]
  \centering
  \begin{tabular}{llllllllll}
    \hline
    Functions & $k_{max}$ & $N$ & $N_{thr}$ & $\lambda$ & $\gamma$ & $Q$ & $\sigma_x$ & $\sigma_y$ & $R$\\\hline\hline
    $H_5(x) = f_1(x)+v$ & 100 & 50 & $N/2$ & 5.5 & 0.15 &  1e-8 & 1e-8 & 1e-8 & 10\\
    $H_6(x) = f_1(x)$ & 100 & 50 & $N/2$ & 3 & 0.15 & 1e-8 & 1e-8 & 1e-8 & 0\\
    $H_7(x) = f_4(x)+v$ & 100 & 50 & $N/2$ & 5.5 & 0.15 &  1e-8 & 1e-8 & 1e-8 & 10 \\
    $H_8(x) = f_4(x)$ & 100 & 50 & $N/2$ &3 & 0.15&  1e-8 & 1e-8 & 1e-8 & 0 \\
    \hline
  \end{tabular}
  \caption{Two CEC 2005 example functions' parameter table}\label{tab:example_func_cec}
\end{table*}


In this section, the PFO algorithm is compared with the PSO algorithm. Since both the benchmark functions and the PSO algorithm assume deterministic measurements, a minor modification is required in the PFO framework to ensure a fair comparison. This adjustment pertains to lines 13 and 14 of Algorithm~\ref{alg:pfo}, where the global minimum and its corresponding measurement are computed.  Given the deterministic nature of the measurements, or at least their determinism in the direction of the minima, the global best solution is selected based on a Maximum a Posteriori (MAP) estimate of the likelihood. It should be emphasized that enforcing a deterministic selection of the global minimum fundamentally alters the behavior of PFO and undermines its original probabilistic interpretation.

Furthermore, it is assumed that the algorithm has access to the true optimal solution for evaluating the termination criterion, specifically requiring a terminal error less than $1\text{e}-8$. The parameter settings used for the PSO algorithm are provided in Table~\ref{tab:pso_param}, and the comparative results over 25 Monte Carlo trials are summarized in Table~\ref{tab:compare}. While both algorithms are capable of identifying the optimal solution, the results show that PSO generally outperforms PFO under these deterministic conditions. However, PSO fails to locate the minima for functions belonging to the class-$\mathcal{C}$ category.

\begin{table*}[ht]
  \centering
  \begin{tabular}{llllllll}
    \hline
    Functions & $N$ & $v_{max}$ & $\phi_{P, max}$ & $\phi_{N, max}$ & $\phi_{G, max}$ & $w_{max}$ & $w_{min}$ \\\hline\hline
    $H_6(x)$ & 150 & 2.26 & 0.37 & 3.68 & 7.4 & 0.9 & 0.25\\
    $H_8(x)$ & 150 & 7.18 & 0.32 & 7.0 & 8.05 & 0.9 & 0.15 \\
    \hline
  \end{tabular}
  \caption{PSO parameters}\label{tab:pso_param}
\end{table*}

\begin{table*}[ht]
  \centering
  \begin{tabular}{llllll}
    \hline
    FES & Criteria & $H_6(x)$, PSO & $H_6(x)$, PFO & $H_8(x)$, PSO & $H_8(x)$, PFO  \\\hline\hline
    1e3 & 1st (best)          & 0.0000 &  0.0000 & 0.0000 & 0.0000\\
        & 7th                 & 0.0000 &  0.0007 & 0.0000 & 0.0051\\
        & 13th (median)       & 0.0006 &  0.0023 & 0.0004 & 0.0103\\
        & 19th                & 0.0019 &  0.0091 & 0.0013 & 0.0256\\
        & 25th (worst)        & 0.0053 &  0.0561 & 0.0060 & 0.3606\\
        & mean                & 0.0012 &  0.0096 & 0.0011 & 0.0416\\
        & std                 & 0.0016 &  0.0151 & 0.0016 & 0.0819\\
    1e4 & 1st (best)          & 0.0000 &  0.0000 & 0.0000 & 0.0000\\
        & 7th                 & 0.0000 &  0.0008 & 0.0000 & 0.0055\\
        & 13th (median)       & 0.0000 &  0.0026 & 0.0000 & 0.0095\\
        & 19th                & 0.0000 &  0.0086 & 0.0000 & 0.0270\\
        & 25th (worst)        & 0.0000 &  0.0556 & 0.0002 & 0.3410\\
        & mean                & 0.0000 &  0.0097 & 0.0000 & 0.0409\\
        & std                 & 0.0000 &  0.0149 & 0.0000 & 0.0792\\
    \hline
  \end{tabular}
  \caption{Comparison of PSO and PFO for functions $H_6(x)$ and $H_8(x)$}\label{tab:compare}
\end{table*}

\subsection{Example Set 3: Multivariable Optimization}

This set of examples evaluates the performance and robustness of the proposed PFO algorithm on multivariable functions. Robustness is assessed through 10 Monte Carlo trials. The algorithm parameters, listed in Table~\ref{tab:example_func_multi}, are tuned using random sampling over predefined set.

\begin{figure}[ht]
      \centering
      \subfloat[]{\label{}\includegraphics[width=1\columnwidth]{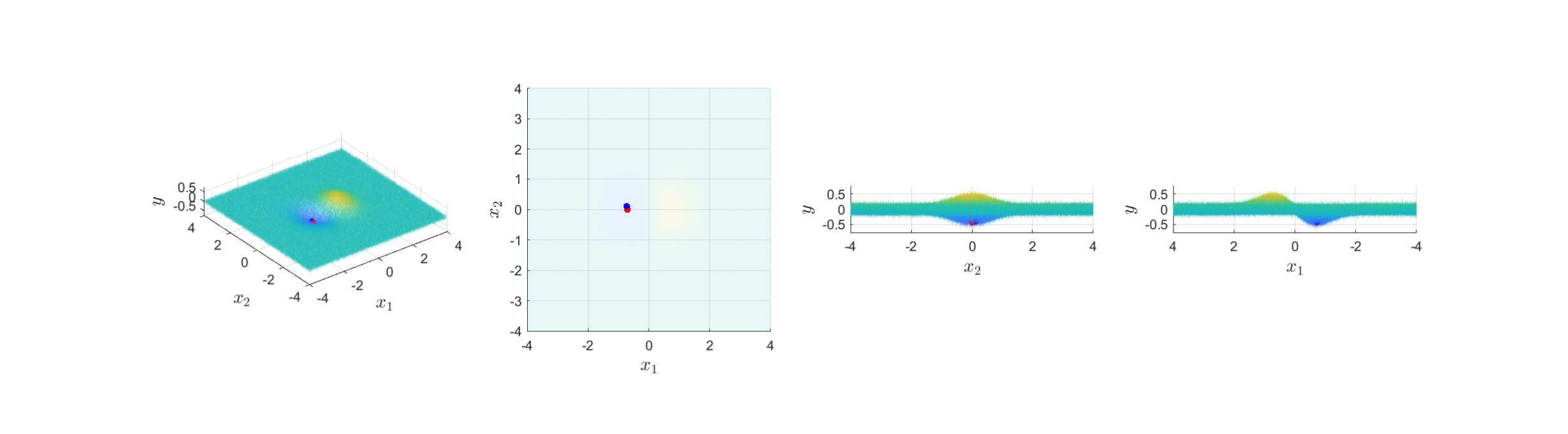}}\newline
      \subfloat[]{\label{}\includegraphics[width=0.6\columnwidth]{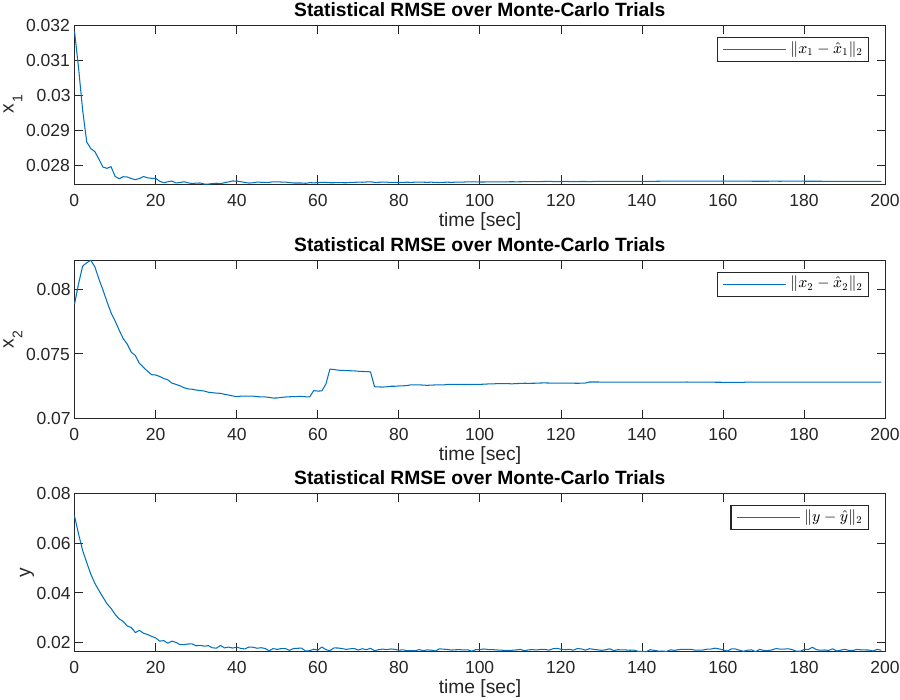}}\\
     \caption{Result of PFO for function $H_9(x)$}\label{fig:h9}
  \end{figure}

\begin{figure}[ht]
  \centering
  \subfloat[]{\label{}\includegraphics[width=1\columnwidth]{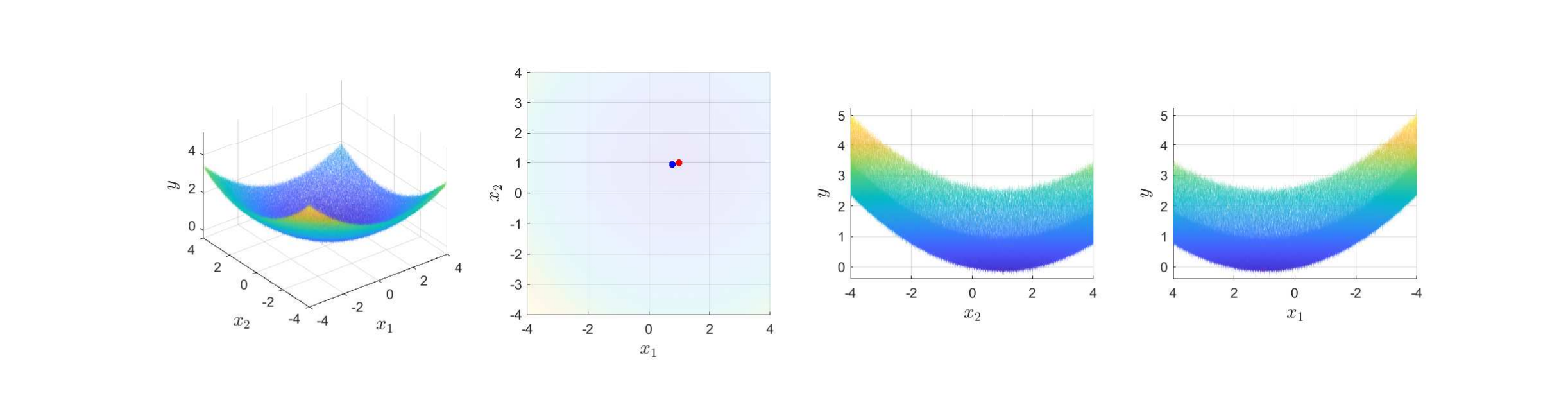}}\newline
  \subfloat[]{\label{}\includegraphics[width=0.6\columnwidth]{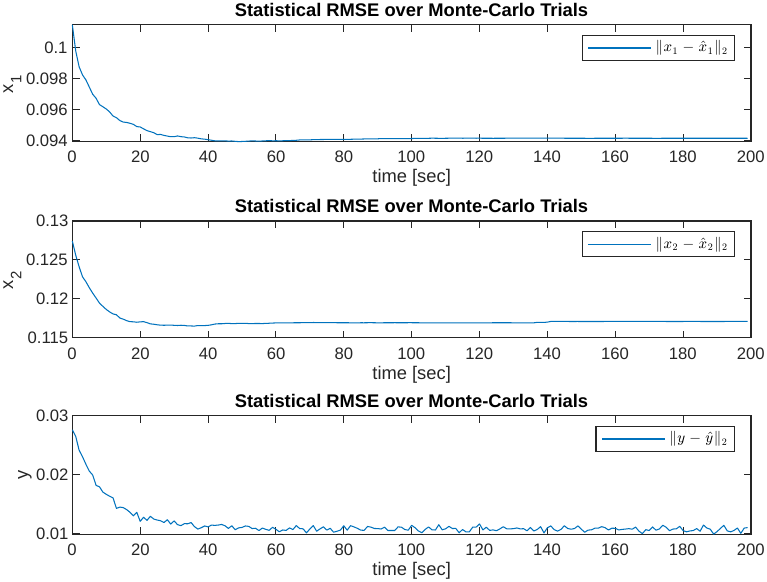}}\\
  \caption{Result of PFO for function $H_{10}(x)$}\label{fig:h10}
\end{figure}

\begin{figure}[ht]
  \centering
  \subfloat[]{\label{}\includegraphics[width=1\columnwidth]{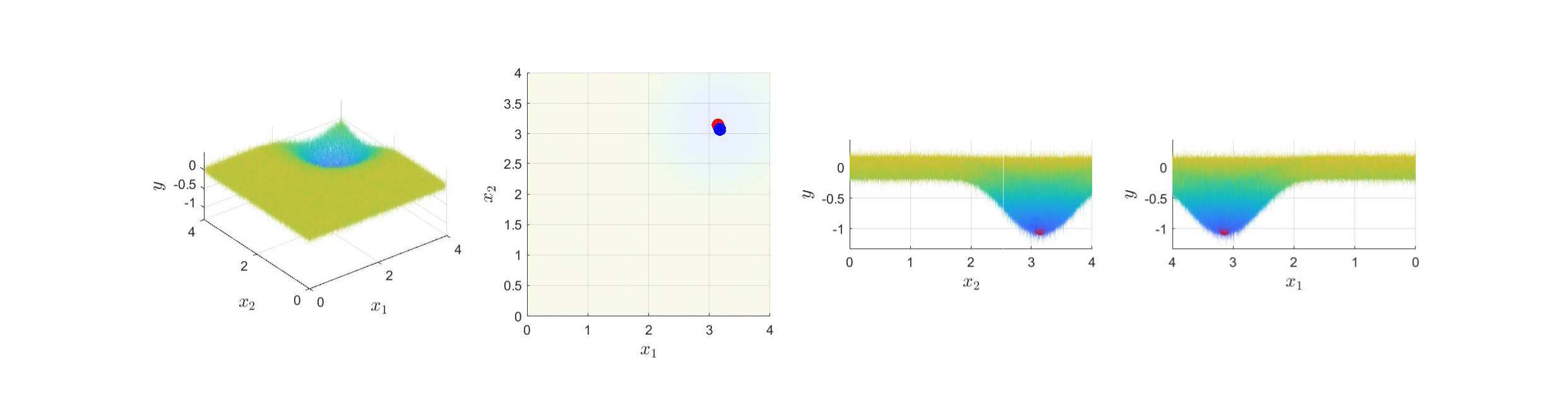}}\newline
  \subfloat[]{\label{}\includegraphics[width=0.6\columnwidth]{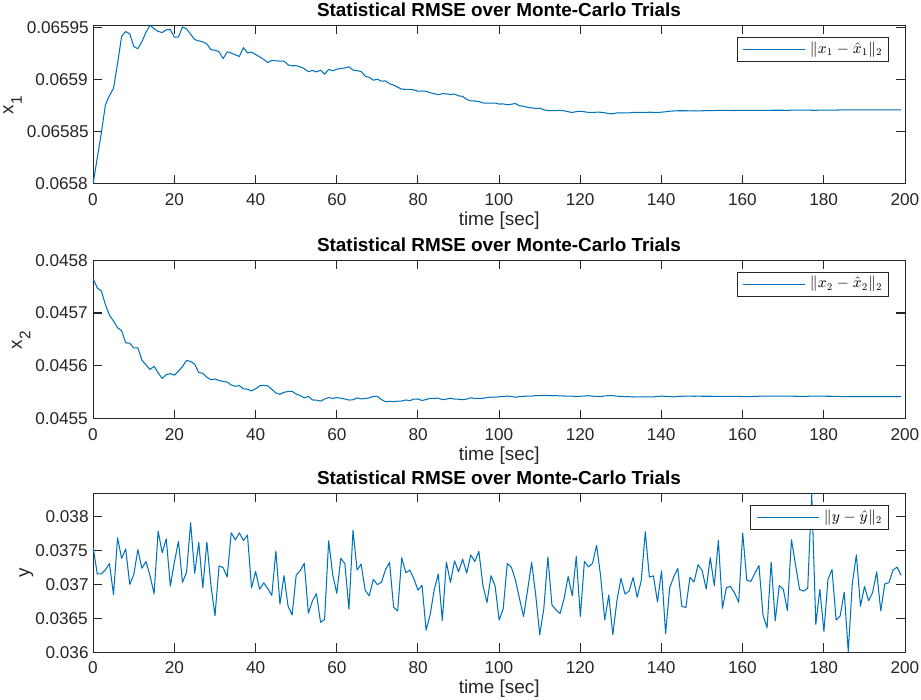}}\\
 \caption{Result of PFO for function $H_{11}(x)$}\label{fig:h11}
\end{figure}

\begin{figure}[ht]
  \centering
  \subfloat[]{\label{}\includegraphics[width=1\columnwidth]{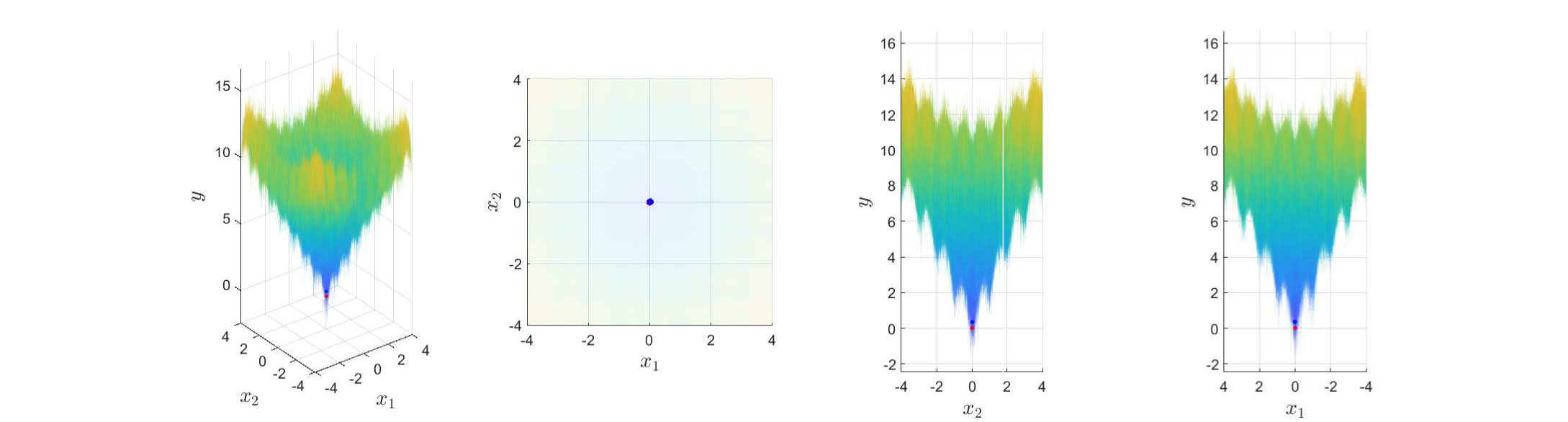}}\newline
  \subfloat[]{\label{}\includegraphics[width=0.6\columnwidth]{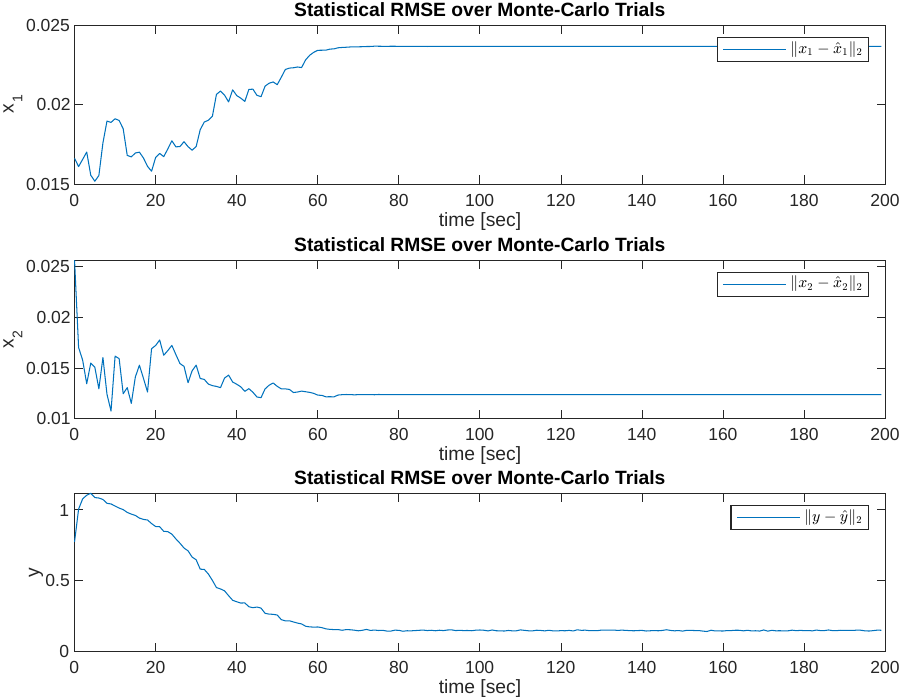}}\\
 \caption{Result of PFO for function $H_{12}(x)$}\label{fig:h12}
\end{figure}

Figures~\ref{fig:h9} to \ref{fig:h12} illustrate the results. Subfigures~(a) show the distribution of the problem data and the best-found solution, while subfigures~(b) depict the statistical RMSE at each step across the Monte Carlo runs. The best-found minimum is marked with a small blue sphere, and the true global minimum is indicated by a small red sphere. The plots demonstrate a high level of confidence in locating the global minimum within the specified uncertainty bounds. Since these functions involve two optimization variables, the weighting matrix is defined as $Q = \text{diag}([q_1, q_2])$. The functions \( H_9(x) \) through \( H_{12}(x) \) are defined by the following equations.

\begin{align}
   & H_9(x) = x_1e^{-x_1^2-x_2^2}+v\nonumber\\
   & H_{10}(x) = -0.1((x_1-1)^2+(x_2-1)^2)\nonumber\\
   & H_{11}(x) = -\cos(x_1)\cos(x_2)e^{-((x_1 - \pi)^2 + (x_2 - \pi)^2)}+v\nonumber\\
   & H_{12}(x) = -20 e^{-0.2 \sqrt{0.5 (x_1^2 + x_2^2)}} - e^{0.5(\cos(2\pi x_1) + \cos(2\pi x_2))} + 20+v
\end{align}
\begin{table*}[ht]
  \centering
  \begin{tabular}{llllllllll}
    \hline
    Functions & $k_{max}$ & $N$ & $N_{thr}$ & $\lambda$ & $\gamma$ & $[q_1,q_2]$ & $\sigma_x$ & $\sigma_y$ & $R$\\\hline\hline
    $H_9(x)$ & 250 & 200 & $N/2$ & 1 & 0.75 &  [1e-4,1e-2] & 3e-9 & 3e-9 & 0.1\\
    $H_{10}(x)$ & 350 & 200 & $N/2$ & 1.5 & 0.75 &  [1e-3,1e-3] & 3e-9 & 3e-9 & 0.1\\
    $H_{11}(x)$ & 250 & 200 & $N/2$ & 1 & 1 &  [1e-4,1e-4] & 3e-9 & 3e-9 & 0.1 \\
    $H_{12}(x)$ & 400 & 200 & $N/2$ & 1 & 0.4&  [3e-3,1e-4] & 3e-9 & 3e-9 & 1 \\
    \hline
  \end{tabular}
  \caption{Multivariable optimization function set's parameter table}\label{tab:example_func_multi}
\end{table*}

  \subsection{Parameter analysis}
\begin{figure}[ht]
  \centering
  \subfloat[]{\label{}\includegraphics[width=0.3\columnwidth]{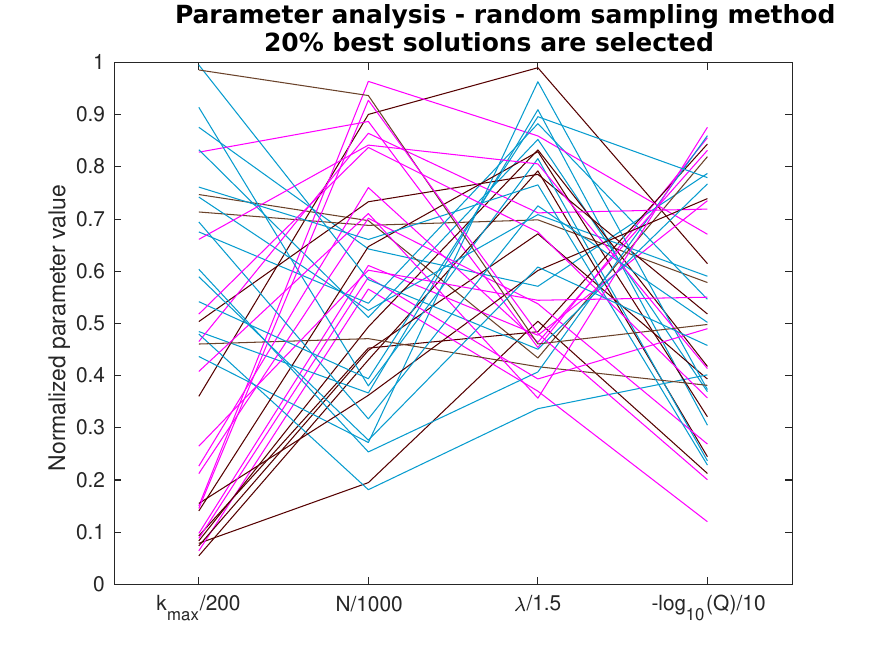}}
  \subfloat[]{\label{}\includegraphics[width=0.3\columnwidth]{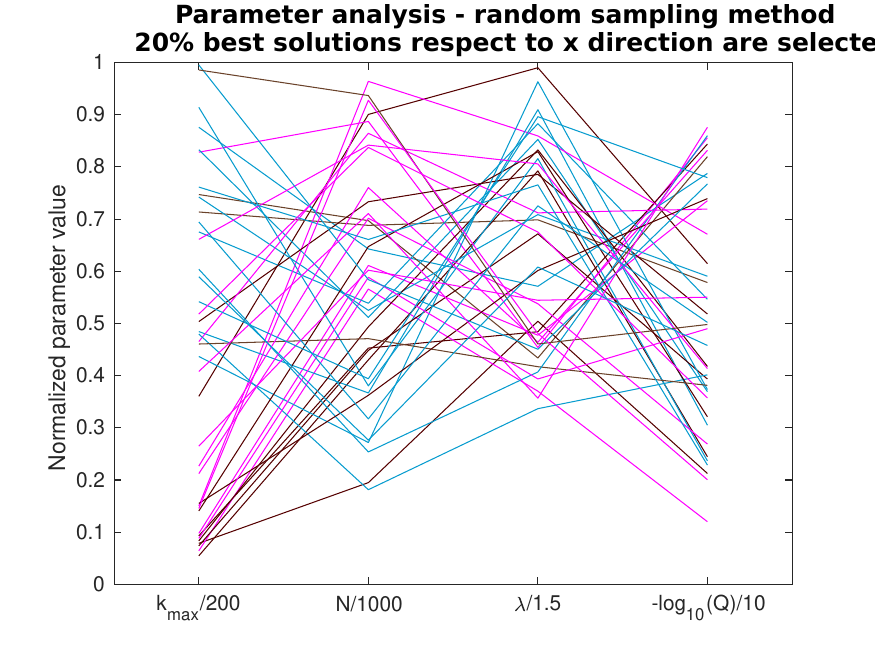}}
  \subfloat[]{\label{}\includegraphics[width=0.3\columnwidth]{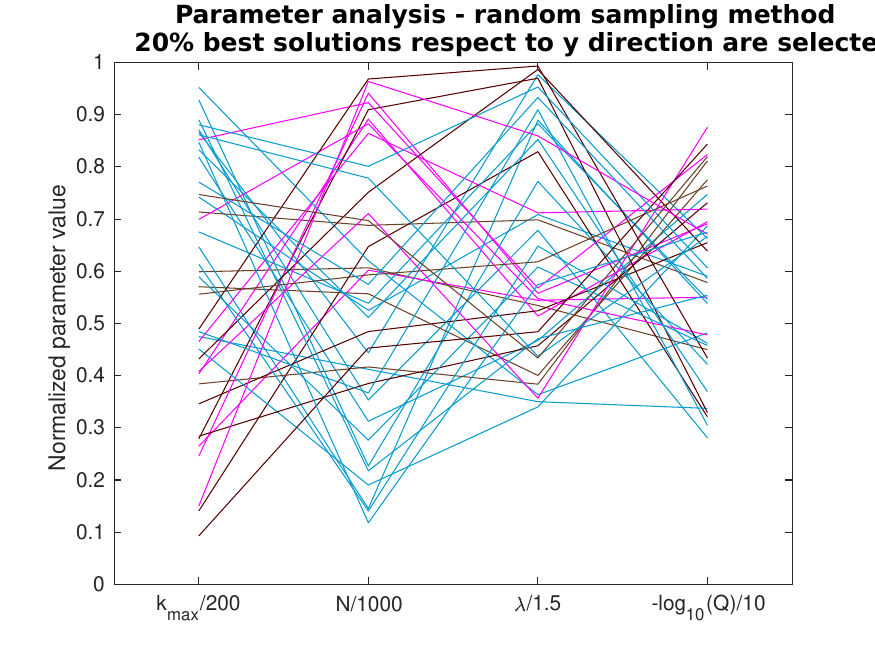}}\\
 \caption{Random sampling results}\label{fig:random_sampling}
\end{figure}

In order to analyze the effect of parameter selection on the optimization error, a random sampling method is utilized \cite{franken2009visual}. Results of random sampling for 200 samples are depicted in Figure~\ref{fig:random_sampling}. The results suggest the following conjectures in parameter selection for single variable optimizations. Trials shown, the results can be extended for relatively higher dimensions, as well.
\begin{conjecture}\label{pos:p1}
It is likely to reach a better solution with \emph{lower} maximum number of iterations and UT scaling factor ($k_{max}$ and $\lambda$) when the number of particles and state transition covariance ($N$ and $Q$) are \emph{high}. Results suggest a correlation between $k_{max}$ and $N$, and between $\lambda$ and $Q$.
\end{conjecture}
\begin{conjecture}
It is likely to reach a better solution with \emph{higher} maximum number of iterations and UT scaling factor ($k_{max}$ and $\lambda$) when the number of particles and state transition covariance ($N$ and $Q$) are \emph{low}.
\end{conjecture}
\begin{conjecture}\label{pos:p3}
  A \emph{moderate} choice for maximum number of iterations and UT scaling factor ($k_{max}$ and $\lambda$) is likely to result in a better solution when the number of particles and state transition covariance ($N$ and $Q$) are \emph{moderate}.
\end{conjecture}
\begin{conjecture}\label{pos:p4}
The patterns given in Conjectures \ref{pos:p1} to \ref{pos:p3} are more important than the parameter values if a good initial guess is considered in the chain of parameters. This leads to low sensitivity to parameter variations.
\end{conjecture}
\begin{conjecture}
Conjectures \ref{pos:p1} to \ref{pos:p4} are the same for optimization errors in each direction (i.e. $x$, $y$ and $(x,y)$).
\end{conjecture}

\section{Conclusion}\label{sec:conc}
It was provisioned from the beginning to adhere theoretical supports to heuristic global optimization algorithm in stochastic optimization problem. The proposed algorithm in this paper under highlighted assumptions nailed this objective. Although, more work is needed to modify the weight update equation based on given information in the transitional prior to maximize this support. The algorithm benefits from a low number of parameters, which are easily tuned based on basic conjectures.  Eventually, three problem sets are attempted by the PFO with promising results and based on Monte-Carlo trials, it has shown robustness. Since the intuition behind the PFO development is based on the uncertainty in the measurements, performance downgrade for deterministic problems is evident. The performance of the PFO is compared with the PSO in noise-free problems. The PSO showed statistically better performance; however, the best-found solution over Monte-Carlos was the same. Indeed, the PSO fails to optimize stochastic class-$\mathcal{C}$ problems.

\bibliographystyle{unsrt}

\bibliography{ref}


\end{document}